\documentclass{TAC}

\usepackage[usenames,dvipsnames]{pstricks}

\usepackage{amscd,amssymb,amsmath,graphicx,verbatim,hyperref}
\usepackage[OT1,T1]{fontenc}

\usepackage[all, cmtip]{xy}
\usepackage{pdfpages}

\usepackage{amsmath,amsfonts,mathabx,tikz,tikz-cd,stmaryrd,url,mathtools, mathrsfs, tensor, pigpen}
\usepackage{hyperref}
\hypersetup{
	colorlinks,
	citecolor=black,
	filecolor=black,
	linkcolor=black,
	urlcolor=black
}

\usetikzlibrary{matrix,arrows,decorations.pathmorphing}

\newcommand{\po}{\ar@{}[dr]|{\text{\pigpenfont R}}}
\newcommand{\pb}{\ar@{}[dr]|{\text{\pigpenfont J}}}
\newcommand{\fpb}{\ar@{}[dr]|{\text{fakepb}}}

\newcommand{\la}{\leftarrow}
\newcommand{\ra}{\rightarrow}
\newcommand{\lra}{\longrightarrow}
\newcommand{\Ra}{\Rightarrow}

\newcommand{\ldual}[1]{\mathord{{\let\nolimits\relax\sideset{^\wedge}{}{#1}}}}
\newcommand{\laction}[2]{\mathord{{\let\nolimits\relax\sideset{^{#1}}{}{#2}}}}
\newcommand{\conj}[2]{\mathord{{\let\nolimits\relax\sideset{^{#1}}{}{#2}}}}

\newcommand{\ox}{\otimes}
\newcommand{\x}{\times}
\newcommand{\xra}{\xrightarrow}
\newcommand{\xla}{\xleftarrow}

\newcommand{\xRa}{\xRightarrow}

\def\CA{{\mathscr A}}
\def\CB{{\mathscr B}}
\def\CC{{\mathscr C}}
\def\CD{{\mathscr D}}

\def\CG{{\mathscr G}}
\def\CH{{\mathscr H}}
\def\CI{{\mathscr I}}

\def\CK{{\mathscr K}}

\def\CM{{\mathscr M}}

\def\CV{{\mathscr V}}

\def\CX{{\mathscr X}}

\makeatletter
\newcommand*\bigcdot{\mathpalette\bigcdot@{.5}}
\newcommand*\bigcdot@[2]{\mathbin{\vcenter{\hbox{\scalebox{#2}{$\m@th#1\bullet$}}}}}
\makeatother

\newcommand{\twocong}[2][0.5]{\ar@{}[#2] \save ?(#1)*{\cong}\restore}
\newcommand{\twoeq}[2][0.5]{\ar@{}[#2] \save ?(#1)*{=}\restore}
\newcommand{\ltwocell}[3][0.5]{\ar@{}[#2] \ar@{=>}?(#1)+/r 0.2cm/;?(#1)+/l 0.2cm/^{#3}}
\newcommand{\rtwocell}[3][0.5]{\ar@{}[#2] \ar@{=>}?(#1)+/l 0.2cm/;?(#1)+/r 0.2cm/^{#3}}
\newcommand{\utwocell}[3][0.5]{\ar@{}[#2] \ar@{=>}?(#1)+/d  0.2cm/;?(#1)+/u 0.2cm/_{#3}}
\newcommand{\dtwocell}[3][0.5]{\ar@{}[#2] \ar@{=>}?(#1)+/u  0.2cm/;?(#1)+/d 0.2cm/^{#3}}
\newcommand{\ultwocell}[3][0.5]{\ar@{}[#2] \ar@{=>}?(#1)+/dr  0.2cm/;?(#1)+/ul 0.2cm/^{#3}}
\newcommand{\urtwocell}[3][0.5]{\ar@{}[#2] \ar@{=>}?(#1)+/dl  0.2cm/;?(#1)+/ur 0.2cm/^{#3}}
\newcommand{\dltwocell}[3][0.5]{\ar@{}[#2] \ar@{=>}?(#1)+/ur  0.2cm/;?(#1)+/dl 0.2cm/^{#3}}
\newcommand{\drtwocell}[3][0.5]{\ar@{}[#2] \ar@{=>}?(#1)+/ul  0.2cm/;?(#1)+/dr 0.2cm/^{#3}}

\begin{document}
\title{Monoidal centres and groupoid-graded categories}
\author{Branko Nikoli\'c and Ross Street}
\address{Department of Mathematics, Macquarie University, NSW 2109, Australia}
\eaddress{<branik.mg@gmail.com> <ross.street@mq.edu.au>}
\thanks{The second author gratefully acknowledges the support of Australian Research Council Discovery Grant DP190102432.}

\date{\small{\today}}
\maketitle
\begin{center}
  \itshape Dedicated to Marta Bunge
  \end{center}

\begin{abstract}
We denote the monoidal bicategory of two-sided modules (also called
profunctors, bimodules and distributors) between categories by $\mathrm{Mod}$; 
the tensor product is cartesian product of categories.
For a groupoid $\CG$, we study the monoidal centre $\mathrm{ZPs}(\CG,\mathrm{Mod}^{\mathrm{op}})$ 
of the monoidal bicategory $\mathrm{Ps}(\CG,\mathrm{Mod}^{\mathrm{op}})$
of pseudofunctors and pseudonatural transformations; the tensor product is pointwise.
Alexei Davydov defined the full centre of a monoid in a monoidal category.
We define a higher dimensional version: the full monoidal centre of a monoidale (= pseudomonoid) in a monoidal bicategory $\CM$, and it is a braided monoidale in the monoidal centre $\mathrm{Z}\CM$ of $\CM$.
Each fibration $\pi : \CH \to \CG$ between groupoids provides an example of a full monoidal centre of a monoidale in 
$\mathrm{Ps}(\CG,\mathrm{Mod}^{\mathrm{op}})$.   
For a group $G$, we explain how the $G$-graded categorical structures, as considered by Turaev and Virelizier in order to construct topological invariants, fit into this monoidal bicategory context. 
We see that their structures are monoidales in the monoidal centre of the monoidal bicategory of $k$-linear 
categories on which $G$ acts.
\end{abstract}

  \tableofcontents  
 
  \section*{Introduction}
  
A feature of categories is that they provide a framework for understanding analogies between concepts in different disciplines as examples of the same concept interpreted in different categories. 
Bicategories do the same for categorical concepts, and so on.
An example of this unification process is provided by the second author in Section 4.8 of \cite{Street2012} where the full centre of an algebra, as defined and constructed by Davydov \cite{Davydov2010}, is shown to be a monoidal centre in the sense of \cite{81} in the monoidal bicategory of pointed categories. 
 
 Another aspect of higher categories is the microcosm 
 principle (so named by Baez and Dolan; see footnote 2 on page 100 of \cite{60}).
 For example, monoids as sets with an associative unital multiplication need no 
 mention of categories in their definition. Yet, once we try to examine the structure
 of monoid in terms of functions and commutative diagrams, we see that we are
 using the monoidal category structure on the category $\mathrm{Set}$ of sets.
 The microcosmic phenomenon is that we need a higher dimensional version
 of monoid (monoidal category) to obtain the general notion of monoid.
 What is needed of the general setting for monoidal category (monoidale) is
 an even higher dimensional notion of monoid, namely, a monoidal bicategory.
 This principle occurs with centres as will appear in Sections~\ref{Roc}
 and \ref{Roic}.  
 
 The bicategory which we call $\mathrm{Mod}$ whose morphisms are two-sided modules between categories was
 (up to biequivalence) considered by Marta Bunge in Chapter III, Section 14 of her PhD thesis \cite{Bunge1966};
 also see \cite{Bunge2011}. Up to duality, the morphisms of $\mathrm{Mod}$ have also been called profunctors 
 \cite{Ben1967}, bimodules \cite{LawMetric} and distributors \cite{BenLectures, BenDistrib}.   
 In Section~\ref{Tbocatsm} we establish our notation for monoidal bicategories with $\mathrm{Mod}$
 as prime example and where comonoidales are the promonoidal categories
 of Brian Day \cite{Day1970}. We begin our main theme by providing examples involving groupoids. 
 We find that there is no more difficulty working with a groupoid $\CG$ than with a group $G$;
 indeed, to some extent the several object case makes the situation clearer than working with one object. 

 In Section~\ref{Roc} we review the centre of a monoidal bicategory. 
 We examine the pointwise-monoidal bicategory $\mathrm{Ps}(\CG,\mathrm{Mod}^{\mathrm{op}})$
 of pseudofunctors from a given groupoid $\CG$ regarded as a locally discrete bicategory. 
 Then, in Section~\ref{Roic}, we generalise a little
 the limit construction of monoidal centre to include internal full centres of 
 Alexei Davydov \cite{Davydov2010} and their higher versions. The centre of a monoidale is obtained
 as a representing object for the centre piece construction (in the terminology of \cite{81}).
 The construction involves descent categories of certain pseudocosimplicial categories which
 arise formally from the universal property of the augmented simplicial category.  
 
 As in the work of \cite{89, 87, Ignacio}, the relationship between centres in $\mathrm{Mod}^{\mathrm{op}}$ and 
 centres of convolution monoidal set-valued functor categories are examined in Section~\ref{CpiModop}.
 This leads naturally in Section~\ref{Tgoaiag} to the study of the groupoid $\CG^{\mathrm{aut}}$ of
 automorphisms in a groupoid $\CG$ since it allows the identification of centres of the 
 cartesian monoidal category $[\CG,\mathrm{Set}]$ and the pointwise-monoidal bicategory
 $\mathrm{Ps}(\CG,\mathrm{Mod}^{\mathrm{op}})$. 
 
 Our main examples arise when we are given a fibration $\pi : \CH\to \CG$ between groupoids
 (such as a surjective group morphism). There is a pseudofunctor $\mathbb{H}$ whose Grothendieck
 construction is $\pi$ and which possesses a natural monoidale structure as an object of
 $\mathrm{Ps}(\CG,\mathrm{Mod}^{\mathrm{op}})$.  
 We discuss the convolution 
 monoidal structure on $\mathrm{Ps}(\CG^{\mathrm{aut}},\mathrm{Mod}^{\mathrm{op}})$
 in Section~\ref{Micb}. 
 There is a pseudofunctor $\mathbb{H}^{\mathrm{aut}}$ whose Grothendieck
 construction is the fibration $\pi^{\mathrm{aut}} : \CH^{\mathrm{aut}}\to \CG^{\mathrm{aut}}$ 
 and which possesses a natural braided monoidale structure as an object of
 $\mathrm{Ps}(\CG^{\mathrm{aut}},\mathrm{Mod}^{\mathrm{op}})$. 
   
 In Sections~\ref{tTVs} and \ref{Ihbab}, we show how structures and constructions defined by 
 Turaev and Virelizier \cite{TurVir2013} 
 (also see \cite{Tur2000, KirPri2008, BarSch2011}) for $G$-graded categories, 
 where $G$ is a group, are precisely the usual monoidal structures taken
 in the monoidal bicategories under discussion. 
 To some extent, we replace the monoidal category of $k$-modules by the monoidal bicategory 
 $\mathrm{Mod}^{\mathrm{op}}$ whose morphisms are two-sided modules on which categories act.
  
 There are several benefits gained by working at the bicategorical level. 
 While $k$-modules only have monoidal duals
 if they are finitely generated and projective, all objects of $\mathrm{Mod}^{\mathrm{op}}$
 have duals. While $k$-modules only have finite direct sums, all small direct sums exist in $\mathrm{Mod}^{\mathrm{op}}$. While the monoidal centre of the category of all $k$-linear representations of a
 group $G$ is not tortile \cite{ShumPhD} (also called ``ribbon''), the monoidal centre of the bicategory of $\mathrm{Mod}^{\mathrm{op}}$-representations of $G$ is tortile. 
 At this stage we are not saying anything about the main theorem of \cite{TurVir2013} which is a $G$-graded
 version of a theorem of M\"uger \cite{Mueger}. 
 However, note that the invertibility of the $S$-matrix in the $G$-graded case uses M\"uger's result.   
 
 Moreover, a higher version of the Davydov full centre occurs in the key example in \cite{TurVir2013}.
 Indeed, our concluding Section~\ref{Afcfagf} shows how $\mathbb{H}^{\mathrm{aut}}$ is the full centre of $\mathbb{H}$.    
   
 \section{The bicategory of categories and two-sided modules}\label{Tbocatsm}

Bicategories as defined by Jean B\'enabou \cite{Ben1967} are ``monoidal categories with several objects'' 
in the sense that additive categories are ``rings with several objects'' \cite{BMitchell}. 
For bicategories $\CA$ and $\CB$, we write $\mathrm{Ps}(\CA,\CB)$ for the bicategory of pseudofunctors $\CA\to\CB$, 
pseudonatural transformations, and modifications (in the terminology of Kelly-Street \cite{Kelly1974}). 
In the spirit of \cite{99} notice that, if $\CA = \CG$ is a groupoid, lax natural transformations between
pseudofunctors $\CA\to\CB$ are automatically pseudonatural. 
We will use monoidal bicategory terminology from Day-Street \cite{60} except that we now use ``monoidale'' in 
preference to ``pseudomonoid''.

A good example of an autonomous symmetric monoidal bicategory is $\mathrm{Mod}$. 
The objects are categories.
The homcategories are defined by 
$$\mathrm{Mod}(\CA,\CB) = [\CB^{\mathrm{op}}\times \CA,\mathrm{Set}] \ ; $$
objects of these homs are called {\em modules}.
Composition $$\mathrm{Mod}(\CB,\CC)\times \mathrm{Mod}(\CA,\CB)\xra{\circ}\mathrm{Mod}(\CA,\CC)$$
is defined by  
$$(N\circ M)(C,A)=\int^B{M(B,A)\times N(C,B)}  \ .$$
The identity module $\CA \to \CA$ is the hom functor of $\CA$. 

Tensor product is finite cartesian product of categories; it is not the product in $\mathrm{Mod}$.
Coproduct in $\mathrm{Mod}$ is coproduct of categories; it is also product in $\mathrm{Mod}$.   

We identify each functor $F : \CA\to \CB$ with the module $F_* : \CA\to \CB$ defined by
$F_*(B,A) = \CB(B,FA)$. There is also the module $F^* : \CB\to \CA$ defined by
$F^*(A,B) = \CB(FA,B)$ and providing a right adjoint $F_*\dashv F^*$ for $F_*$ in $\mathrm{Mod}$.   

The dual of $\CA$ in $\mathrm{Mod}$ is the opposite category $\CA^{\mathrm{op}}$.
We have an equivalence of pseudofunctors in the square
 \begin{equation}\label{starmodules}
 \begin{aligned}
\xymatrix{
\mathrm{Cat}^{\mathrm{co}} \ar[d]_{(-)^{\mathrm{op}}}^(0.5){\phantom{aaaaaa}}="1" \ar[rr]^{(-)^*}  && \mathrm{Mod}^{\mathrm{op}} \ar[d]^{(-)^{\mathrm{op}}}_(0.5){\phantom{aaaaaa}}="2" \ar@{=>}"1";"2"^-{\simeq}
\\
\mathrm{Cat} \ar[rr]_-{(-)_*} && \mathrm{Mod} 
}
\end{aligned}
\end{equation} 
where the top $(-)^*$ and bottom $(-)_*$ pseudofunctors are the identity on objects, locally-full, strong monoidal, and coproduct preserving. The right side of \eqref{starmodules} is an equivalence (bijective on objects and an isomorphism
on homcategories).

For each category $\CA$, we write $\mathrm{Q}\CA$ for the splitting idempotent completion of $\CA$ 
(for example, see Chapter 2, Exercises B of \cite{Freyd1966}) with inclusion $\mathrm{N}_{\CA} :\CA\to \mathrm{Q}\CA$. 
It is easy to see that $(\mathrm{Q}\CA)^{\mathrm{op}}\simeq \mathrm{Q}(\CA^{\mathrm{op}})$ and $\mathrm{Q}\CA\times \mathrm{Q}\CB \simeq \mathrm{Q}(\CA\times \CB)$.  
Since idempotent splittings are preserved by all functors, we also have $[(\mathrm{Q}\CA)^{\mathrm{op}},\mathrm{Set}]\simeq [\CA^{\mathrm{op}},\mathrm{Set}]$. 
Recall the non-additive Morita-type Theorem from \cite{Bunge1966} that $[\CA^{\mathrm{op}},\mathrm{Set}]\simeq [\CB^{\mathrm{op}},\mathrm{Set}]$ if and only if $\mathrm{Q}\CA\simeq \mathrm{Q}\CB$. 
Consequently, we have a strong monoidal auto-biequivalence on $\mathrm{Mod}$: 
\begin{eqnarray}\label{Q}
\mathrm{Q} : \mathrm{Mod}\to \mathrm{Mod}
\end{eqnarray}
and a monoidal pseudonatural equivalence $\mathrm{N} : 1_{\mathrm{Mod}} \simeq \mathrm{Q}$. 
From \cite{LawMetric} we also see that $\mathrm{Q}$ has a Cauchy completion property which implies it takes left adjoint modules (``Cauchy modules'') to functors (``convergent modules''). In particular, equivalences in $\mathrm{Mod}$ are
taken to equivalences in $\mathrm{Cat}$.

Recall (for example from \cite{60}) that monoidales in $\mathrm{Mod}^{\mathrm{op}}$ are promonoidal
categories $\CA$ in the sense of Day \cite{Day1970}. The tensor product is a module $P : \CA \to \CA\times \CA$
and the unit is a module $J : \CA \to \mathbf{1}$. The Day convolution monoidal structure on the functor category $[\CA,\mathrm{Set}]$ has unit $J : \CA\to \mathrm{Set}$ and tensor product $F\star G$ defined by
\begin{eqnarray*}
(F\star G)C = \int^{A,B}{P(A,B;C)\times FA\times GB} \ .
\end{eqnarray*}

Suppose $\CV$ is a complete cocomplete closed symmetric monoidal category.
We have a symmetric (weak) monoidal pseudofunctor
\begin{eqnarray}\label{[-,V]}
[-,\CV] : \mathrm{Mod}^{\mathrm{op}} \lra \CV\text{-}\mathrm{CAT}
\end{eqnarray}
taking each category $\CA$ to the $\CV$-enriched functor category $[\CA,\CV]$.
Therefore, each promonoidal category is taken to a convolution monoidal $\CV$-category. 
In the case $\CV = \mathrm{Set}$, \eqref{[-,V]} is the contravariant representable pseudofunctor 
$\mathrm{Mod}(-,\mathbf{1})$.     
\bigskip

Let $\CG$ be a groupoid; that is, a category in which all morphisms are invertible. 
Like every category, $\CG$ has a promonoidal structure for which the convolution tensor-product structure on $[\CG,\mathrm{Set}]$
is pointwise cartesian product. In this case, the tensor unit $J$ is constant at the one-point set and the module $P : \CG \to \CG\times \CG$ is defined by
\begin{eqnarray}\label{promonoidalG}
P(p,q;r) = \CG(p,r) \times \CG(q,r), \ P(a,b;c)(x,y) = (\CG(a,c)x,\CG(b,c)y) = (cxa,cyb).
\end{eqnarray}

As we expound the theory, we will carry along an example based on Example 3.5 of \cite{TurVir2013}
which, in turn, was based on an example in \cite{MNS}.
At this point it is just to show some naturally occurring lax functors and pseudofunctors.
 
\begin{example}\label{ex1}
A groupoid morphism $\pi : \CH\to \CG$ is a functor.
According to B\'enabou \cite{BenLectures} (see \cite{Street2001}), any functor over a category 
$\CG$ corresponds to a normal lax functor $\mathbb{H} : \CG \to \mathrm{Mod}^{\mathrm{op}}$. 
The functor $\pi$ is an iso-fibration if and only if it is a fibration if and only if it is an opfibration if and only if it is Giraud-Conduch\'e \cite{Giraud64, Cond}. In the one-object case, it is the same as saying the group morphism $\pi : H\to G$ is surjective. In the general case, it means that, for all morphisms $p \xra{g} \pi(t)$ in $\CG$, there exists a morphism
$s\xra{h}t$ in $\CH$ such that $\pi(s)=p$ and $\pi(h)=g$. 
(All morphisms in $\CH$ are cartesian and cocartesian for $\pi$.)
In this case, the normal lax functor is actually a pseudofunctor.
Let us now describe it. For $p\in \CG$, the category $\mathbb{H}p$ is the fibre $\pi^{-1}(p)$ of the functor $\pi$ over $p$:
the objects are those $s\in \CH$ with $\pi(s)=p$ and the morphisms those $x :s\to s'$ in $\CH$ with $\pi(x)=1_p$.
For $a : p \to q$ in $\CG$, the module $\mathbb{H}a : \mathbb{H}q \to \mathbb{H}p$ is defined by
\begin{eqnarray*}\label{bfH}
(\mathbb{H}a)(s,t) = \{x \in \CH(s,t) : \pi(x) = a\} \ .
\end{eqnarray*}
In particular, $\mathbb{H}1_p$ is the identity module of $\mathbb{H}p$.
For a composable pair $p\xra{a}q \xra{b}r$ in $\CG$, the composite module $\mathbb{H}a \circ \mathbb{H}b$ is defined by
\begin{eqnarray*}
(\mathbb{H}a \circ \mathbb{H}b)(s,u) = \int^t \{y\in \CH(t,u) : \pi(y) = b\}\times \{x\in \CH(s,t) : \pi(x) = a\} \ .
\end{eqnarray*}
The composition constraint has components the bijections  
      \begin{eqnarray*}
(\mathbb{H}a \circ \mathbb{H}b)(s,u) \lra \mathbb{H}(ba)(s,u) \ , \ [y,x] \mapsto yx \ .
\end{eqnarray*}

A cleavage $\sigma$ for the fibration $\pi : \CH\to \CG$ amounts to a choice of morphism $a^*(t)\xra{\sigma(a)}t$ 
in $\CH$ with $\pi(\sigma(a)) = a$ for each $t\in \CH$ and $a\in\CG(p,\pi(t))$.
The pseudofunctor $\mathbb{H} : \CG \to \mathrm{Cat}^{\mathrm{op}}$ corresponding to this cloven fibration is defined
by taking $\mathbb{H}p = \pi^{-1}(p)$ as before, taking $\mathbb{H}a = a^*$ where $a^*(t\xra{y}t')= (a^*(t)\xra{\sigma(a)}t\xra{y}t'\xra{\sigma^{-1}(a)}a^*(t'))$, while the component at $u\in \mathbb{H}r$ of the invertible composition constraint $(\mathbb{H}a)(\mathbb{H}b)\Rightarrow \mathbb{H}(ba)$ is equal to $\sigma(ba)^{-1}\sigma(b)\sigma(a)\in (\mathbb{H}p)(a^*b^*(u), (ba)^*(u))$. 
We use the same symbol for the pseudofunctors $\mathbb{H}$ since the first is equivalent to the composite 
of the second with the canonical pseudofunctor $(-)_*$.  
\end{example}  

\begin{remark} Fr\"ohlich-Wall \cite{FW} studied monoidales in the cartesian monoidal 2-category $\mathrm{Cat}/\CG$. 
(For us $\CG$ can be a general groupoid although, for them, $\CG$ has one object, that is, $\CG$ is a group). 
They were mainly interested in objects $\pi : \CC\to \CG$ for which $\CC$ is a categorical group (what they call ``group-like'')
and $\pi$ is Giraud-Conduch\'e.
So they could equally have looked at the monoidal full sub-2-category $\mathrm{Cat}/_{\mathrm{GC}}\CG$ of $\mathrm{Cat}/\CG$ whose objects are Giraud-Conduch\'e functors.
There is a strong monoidal pseudofunctor $$\mathrm{Cat}/_{\mathrm{GC}}\CG \to \mathrm{Ps}(\CG,\mathrm{Mod}^{\mathrm{op}})$$ 
which takes morphisms to pseudonatural transformations whose components are right adjoints of functors;
for this see \cite{Street2001} and note that taking inverses of morphisms defines a (promonoidal) 
isomorphism $\CG^{\mathrm{op}}\cong \CG$ and we have the monoidal equivalence $(-)^{\mathrm{op}} : \mathrm{Mod}^{\mathrm{op}}\simeq \mathrm{Mod}$
of \eqref{starmodules}.    
\end{remark}

\section{Review of centres}\label{Roc} 

The centre of a monoidal category appeared in \cite{38} and was reported by S. Majid \cite{Majid} as known to
V. Drinfel'd along with its connection to the Drinfel'd double \cite{Drinfeld} of a Hopf algebra. The centre\footnote{We like the fact that $\mathrm{Z}$ is not only the first letter of the German word for centre but also for braid.} $\mathrm{Z}\CV$ of a
monoidal category $\CV$ is a braided monoidal category in the sense of \cite{Joyal1993}. 

The centre $\mathrm{Z}\CM$ of a monoidal bicategory $\CM$ appeared in \cite{BN1996} with corrections in \cite{Crans1998} and applications to topological quantum field theory in \cite{KTZ}.
The centre of a monoidal bicategory is a braided monoidal bicategory in the sense of \cite{89}.

Taking a monoidal bicategory $\CM$ to be a one-object tricategory in the sense of \cite{GPS}, McCrudden \cite{Paddy} explicitly defined braidings on $\CM$ in his Appendix A. However, we will take the approach of 
\cite{BN1996, Crans1998, 89} and work as if our monoidal bicategory $\CM$ were a Gray monoid
(also called a semistrict monoidal bicategory). Justification lies in the coherence theorem of \cite{GPS}.
In particular, we write as if $\CM$ were a (strict) 2-category, and associativity and unit constraints 
(along with their equivalence adjoints, units and counits) were identities.
One of the two choices will be assumed made to obtain the tensor product $-\ox - :\CM\times \CM \to \CM$ 
as a pseudofunctor. In some diagrams we will even delete the tensor symbol $\ox$ between its arguments.
The tensor unit is denoted by $I$.     

An object of $\mathrm{Z}\CM$ is a triplet $\underline{A} = (A,u,\zeta)$ where $A$ is an object of $\CM$, where $u : A\ox -\to -\ox A$
is a pseudonatural equivalence with $u_I = 1_A$, and $\zeta$ is an invertible modification 
\begin{equation}
\begin{aligned}
\xymatrix{
A\ox -\ox \sim \ar[rr]^-{u} \ar[rd]_{u\ox 1} && -\ox \sim \ox A  \\
& -\ox A\ox \sim \ar[ru]_{1\ox u} \dtwocell[0.55]{u}{\zeta}
}
\end{aligned}
\end{equation}
 with $\zeta_{I,X}= 1_{u_X} = \zeta_{X,I}$ subject to the 2-cocycle condition \eqref{3simplex}.  
\begin{equation}\label{3simplex}
 \begin{aligned}
 \xymatrix{
AXYZ\ar[rrdd]^{u_{XY}1_Z}\ar[rr]^{u_{XYZ}} \ar[dd]_{u_X1_Y1_Z} && \dtwocell[0.45]{lld}{\zeta_{XY,Z}}  \dtwocell[0.7]{lldd}{\zeta_{X,Y}1_Z} XYZA \\
    \\
 XAYZ \ar[rr]_{1_Xu_Y1_Z} &&XYAZ\ar[uu]_{1_X1_Yu_Z}
} 
\qquad 
\xymatrix{ \\
= \\
} 
\qquad
\xymatrix{
AXYZ\ar[rr]^{u_{XYZ}} \ar[dd]_{u_X1_Y1_Z} && \dtwocell[0.75]{lld}{\zeta_{X,YZ}}  \dtwocell[0.75]{ldd}{1_X\zeta_{Y,Z}} XYZA \\
    \\
 XAYZ\ar[rruu]_{1_Xu_{YZ}} \ar[rr]_{1_Xu_Y1_Z} &&XYAZ\ar[uu]_{1_X1_Yu_Z}
}
\end{aligned}
\end{equation}

A morphism $(f,\phi) : (A,u,\zeta)\to (B,v,\xi)$ in $\mathrm{Z}\CM$ consists of a morphism $f : A\to B$ in $\CM$
and an invertible modification
\begin{equation}
\begin{aligned}
\xymatrix{
A\ox -\ar[rr]^-{u} \ar[d]_{f\ox 1} && -\ox A \ar[d]^{1\ox f}  \\
 B\ox - \ar[rr]_{v} && -\ox B \dtwocell[0.5]{llu}{\phi}
}
\end{aligned}
\end{equation}
with $\phi_I = 1_f$ subject to the condition \eqref{tricylinder}.
\begin{equation}\label{tricylinder}
 \begin{aligned}
 \xymatrix{
AXY\ar[rr]^{u_{XY}}\ar[rd]^{u_{X}1_Y} \ar[dd]^{f1_X1_Y} &  & \dtwocell[0.5]{lld}{\zeta_{X,Y}}  XYA \ar[dd]^{1_X1_Yf} \\
& XAY \ar[ru]_{1_Xu_{Y}} \ar[dd]^{1_X f 1_Y} &
    \\ BXY \dtwocell[0.5]{ru}{\phi_{X}1_Y} \ar[rd]_{v_X1_Y} &  &XYB \dtwocell[0.5]{lu}{1_X\phi_{Y}}\\
 & XBY \ar[ru]_{1_Xv_Y} &
} 
\qquad 
\xymatrix{ \\
= \\
} 
\qquad
\xymatrix{
AXY\ar[rr]^{u_{XY}} \ar[dd]^{f1_X1_Y} &  &   XYA \ar[dd]^{1_X1_Yf} \\
&  &
    \\ BXY \ar[rr]^{v_{XY}} \dtwocell[0.5]{rruu}{\phi_{XY}} \ar[rd]_{v_X1_Y} &  &XYB \\
 & XBY \dtwocell[0.55]{u}{\xi_{X,Y}} \ar[ru]_{1_Xv_Y} &
} 
\end{aligned}
\end{equation} 
A 2-morphism $\sigma : (f,\phi)\Rightarrow (g,\psi) : \underline{A}\to \underline{B}$ in $\mathrm{Z}\CM$ is a 2-morphism
$\sigma : f\Rightarrow g$ in $\CM$ such that \eqref{2morcentre} commutes.
\begin{equation}\label{2morcentre}
\begin{aligned}
\xymatrix{
(1_X\ox f)\circ u_X \ar[rr]^-{\phi_X} \ar[d]_-{(1_X\ox \sigma)\circ u_X} && v_X\circ (f\ox 1_X) \ar[d]^-{v_X\circ (\sigma\ox 1_X)} \\
(1_X\ox g)\circ u_X \ar[rr]_-{\psi_X} && v_X\circ (g\ox 1_X)}
\end{aligned}
\end{equation}

The tensor product of $\mathrm{Z}\CM$ is defined on objects by
\begin{eqnarray*}
\underline{A}\ox \underline{A'} = (A\ox A', A\ox A'\ox - \xra{1_A\ox u'}A\ox -\ox A'\xra{u\ox 1_{A'}}-\ox A\ox A', \theta)
\end{eqnarray*}
where $\theta_{X,Y}$ is the composite of the 2-morphism
\begin{eqnarray*}
u_{XY}1_{A'}\circ 1_A u'_{XY}\xRightarrow{(\zeta_{X,Y}\ox 1_{A'})\circ (1_A\ox\zeta'_{X,Y})}
1_Xu_Y1_{A'}\circ u_X1_Y1_{A'}\circ 1_A1_Xu'_Y\circ 1_Au'_X1_Y
\end{eqnarray*}
with the canonical isomorphism between the codomain and the morphism
$$1_Xu_Y1_{A'}\circ 1_X1_Au'_{Y}\circ u_X1_{A'}1_Y\circ 1_Au'_X1_Y \ .$$
On morphisms the tensor product is defined by $(f,\phi)\ox (f',\phi') = (f\ox f', \omega)$
where  
\begin{eqnarray*}
1_X f f'\circ u_X 1_{A'}\circ 1_A u'_X\xRightarrow{\omega_X} v_X 1_{A'}\circ 1_A v'_X\circ f f' 1_X
\end{eqnarray*}
is canonically isomorphic to 
\begin{eqnarray*}
1_X f 1_{A'}\circ u_X 1_{A'}\circ 1_A1_X f'\circ 1_A u'_X\xRightarrow{(\phi_X\ox 1_{A'})\circ (1_A\ox\phi')} v_X 1_{A'}\circ f1_A1_{A'}\circ 1_A v'_X\circ 1_A f' 1_X \ .
\end{eqnarray*}

The braiding for $\mathrm{Z}\CM$ involves a pseudonatural transformation with object component
\begin{eqnarray*}
c_{\underline{A}, \underline{A'}} : \underline{A}\ox \underline{A'}\to \underline{A'}\ox \underline{A}
\end{eqnarray*}
made up of the morphism $A\ox A'\xra{u_{A'}} A'\ox A$ and the pasted invertible 2-morphisms
\begin{eqnarray*}
\begin{aligned}
\xymatrix{
AA'X \ar[dd]_{u_{A'}1_X}\ar[rrdd]^{u_{A'X}}\ar[rr]^{1_Au'_{X}} && AXA' \ar[rrdd]^{u_{XA'}} \dtwocell[0.5]{dd}{u_{1_Au'_X}} \ar[rr]^{u_{X}1_{A'}} && XAA' \ar[dd]^{1_Xu_{A'}} \dtwocell[0.6]{ld}{\zeta^{-1}_{X,A'}}\\
&& && \\
A'AX \ar[rr]_{1_{A'}u_{X}} \dtwocell[0.6]{ru}{\zeta_{A',X}} && A'XA \ar[rr]_{u'_X1_{A}} && XA'A \ ,
}
\end{aligned}
\end{eqnarray*}
and with morphism component the pasted composite 
\begin{equation*}
\xymatrix{
AA' \dtwocell[0.5]{rrd}{\phi_{A'}} \ar[rr]^-{u_{A'}} \ar[d]_-{f1_{A'}} && A'A \ar[d]^-{1_{A'}f} \\
BA' \dtwocell[0.5]{rrd}{v_{f'}} \ar[rr]^-{v_{A'}} \ar[d]_-{1_Bf'} && A'B \ar[d]^-{f'1_B} \\
BB' \ar[rr]_-{v_{B'}} && B'B  \ .\\
}
\end{equation*}
The braiding also involves two invertible 2-cells 
$$\rho_{\underline{A} |  \underline{A'}, \underline{A''}} : c_{\underline{A},\underline{A'}\ox \underline{A''}}\Rightarrow (1_{\underline{A'}}\ox c_{\underline{A}, \underline{A''}})\circ (c_{\underline{A}, \underline{A'}}\ox 1_{\underline{A''}})$$ and  $$\rho_{\underline{A}, \underline{A'} | \underline{A''}} : c_{\underline{A}\ox\underline{A'}, \underline{A''}}\Rightarrow (c_{\underline{A}, \underline{A''}}\ox 1_{\underline{A'}})\circ (1_{\underline{A}}\ox c_{\underline{A'}, \underline{A''}}) \ ; $$
these are respectively provided by the following two diagrams.
\begin{equation*}
\begin{aligned}
\xymatrix{
AA'A'' \ar[rr]^-{u_{A'A''}} \ar[rd]_{u_{A'}\ox 1} && A'A''A  \\
& A'AA'' \ar[ru]_{1\ox u_{A''}} \dtwocell[0.55]{u}{\zeta_{A',A''}}
}
\quad
\xymatrix{
AA'A'' \ar[r]^-{1_A\ox u'_{A''}} \ar[rd]_{1_A\ox u'_{A''}} & AA''A' \ar[r]^-{u_{A''}\ox 1_{A'}} & A''AA'  \\
& AA''A' \ar[ru]_{u_{A''}\ox 1_{A'}} \dtwocell[0.55]{u}{=}
}
\end{aligned}
\end{equation*}

If $\CM$ is a braided monoidal bicategory then there is a canonical braided strong monoidal pseudofunctor
\begin{eqnarray}\label{MinZM}
\CM\lra \mathrm{Z}\CM \ , \ X\mapsto \underline{X} : = (X, c_{X,-}, \rho_{X | -,-}) \ , \ f \mapsto \underline{f} : = (f, c_{f,-}) \ , \ \sigma \mapsto \sigma
\end{eqnarray}
which is locally full.

Of course, if $\CV$ is a monoidal category, it can be regarded as a monoidal bicategory with only identity 2-morphisms, and $\mathrm{Z}\CV$ simplifies to the monoidal centre of $\CV$ as in \cite{38}. 
If $A$ is a monoid (in $\mathrm{Set}$), it can be regarded as a monoidal category with only identity morphisms,
and $\mathrm{Z}A = \{a\in A : ax = xa \ \forall \ x\in A\} \cong \tensor*[^A]{\mathrm{Set}}{^A}(A,A)$ where
$\tensor*[^A]{\mathrm{Set}}{^A}$ is the category Eilenberg-Moore algebras for the monad $A\ox - \ox A$
on $\mathrm{Set}$ and $A$ acts on itself on both sides by its own multiplication.

\begin{example}\label{centreofptwise}
For our groupoid $\CG$, we have the monoidal bicategory $\mathrm{Ps}(\CG,\mathrm{Mod}^{\mathrm{op}})$
where the tensor product is pointwise the tensor of $\mathrm{Mod}^{\mathrm{op}}$; this tensor product is of course Day convolution with the promonoidal structure \eqref{promonoidalG} on $\CG$. Up to equivalence, the objects of $\mathrm{Z}\mathrm{Ps}(\CG,\mathrm{Mod}^{\mathrm{op}})$ can be simplified somewhat; a lower dimensional version appears as Proposition 4.3 of \cite{87} and as Theorem 8.6 of \cite{89}.
To see how this works, take such an object $(F,u, \zeta)$. Making use of the biequivalence $\mathrm{Q}$ of \eqref{Q}, we can assume that $F$ is a pseudofunctor from $\CG$ to $\mathrm{Cat}^{\mathrm{op}}$ and that the equivalences $u_K$ have functors as components. 
By the bicategorical Yoneda lemma \cite{14} and the fact
that $(-)_* : \mathrm{Cat} \to \mathrm{Mod}$ preserves bicategorical colimits, the pseudonatural family of functors 
$u_K : K \times F \to F\times K$ is determined by restricting $K$ to representables. 
The functors $u_{\CG(r,-)}p : \CG(r,p)\times Fp  \to Fp\times \CG(r,p)$ correspond to functions
$\CG(r,p)\to [Fp, Fp\times \CG(r,p)]$ which are, in particular, pseudonatural in $r$. 
By the bicategorical Yoneda lemma, these functions correspond to functors $\bar{\delta}_p : Fp\to Fp\times \CG(p,p)$.
The extra structure and axioms on $u$ are equivalent to giving an isomorphism between the identity of $Fp$ and
 the composite of $\bar{\delta}_p$ with the first projection, so that $\bar{\delta}_p$ is determined up to isomorphism
 by functors $\delta_p : Fp\to \CG(p,p)$, and that these $\delta_p$ are pseudonatural in $p\in \CG$. This last means that
 diagram \eqref{dinatdelta} commutes for all $g\in\CG(p,q)$, that the unit constraint $1_{Fp}\cong F1_p$ is identified by $\delta_p$, and
 that the composition constraint $Fg\circ Fh \cong F(hg)$ is identified by $\CG(1_r,gh)\circ \delta_r$.  
 \begin{eqnarray}\label{dinatdelta}
 \begin{aligned}
\xymatrix{
Fq \ar[r]^-{\delta_q} \ar[d]_-{Fg} & \CG(q,q)\ar[r]^-{ \CG(g,1_q)}  & \CG(p,q)  \\
Fp \ar[rr]_-{\delta_p} & & \CG(p,p)\ar[u]_-{\CG(1_p,g)} }
\end{aligned}
\end{eqnarray}
So every object of $\mathrm{Z}\mathrm{Ps}(\CG,\mathrm{Mod}^{\mathrm{op}})$ is equivalent to one of the form
$(F,u, \zeta)$ with $F$ landing in $\mathrm{Cat}^{\mathrm{op}}$, with $u : - \times F \to F\times -$ obtained from a $\delta$ as above via the formula 
\begin{eqnarray}\label{deltavariant}
u_{K,p}(k,x) = (x, K(\delta_p(x))k)
\end{eqnarray}
for $k\in Kp$ and $x\in Fp$, and with $\zeta$ amounting to a rebracketing isomorphism.
\end{example}

\section{Review of internal centres}\label{Roic}

Now we shall define the full centre of a monoidale in a monoidal bicategory and the centre of a monoidale in a braided monoidal bicategory as birepresenting objects, with somewhat more detail than \cite{81}. Also see \cite{Ignacio}.

Let $\Delta$ denote the algebraist's simplicial category: the objects are the ordinals $\underline{n} = \{0, 1, \dots, n-1\}$
(including the empty ordinal $\underline{0}$) and the morphisms are order-preserving functions. Ordinal sum
provides a monoidal structure on $\Delta$.
Let $\Delta_{\bot, \top}$ denote the subcategory of $\Delta$ consisting of the non-empty ordinals and functions
which preserve first and last elements as well as order. There is a canonical isomorphism of categories
$$\Delta^{\mathrm{op}}\cong \Delta_{\bot, \top}$$
taking $\underline{n}$ to $\underline{1+n}$ and $\xi : \underline{n} \to \underline{m}$ in $\Delta$ to the right adjoint  
of the functor $\underline{1+} \xi : \underline{1+n} \to \underline{1+m}$.

Let $A$ be a monoidale (= pseudomonoid) in a monoidal bicategory $\CM$: 
$$A = \lgroup A,A\ox A\xra{P}A,I\xra{J}A, P(P\ox 1_A)\xRa{\Phi}P(1_A\ox P), P(J\ox 1_A)\xRa{\lambda}1_A\xRa{\rho} P(1_A\ox J)\rgroup \ .$$ 
 A monoidal pseudofunctor 
$A^* : \Delta \to \CM$ is defined by $A^*\underline{n} = A^{\ox n}$ and
$$A^*(\underline{0}\ra\underline{1}\la \underline{2}) = (I\xra{J} A \xla{P} A\ox A) \ .$$ 

Write $\tensor*[^A]{\CM}{^A}$ for the bicategory of pseudo-algebras for the pseudomonad $A\ox - \ox A$
on $\CM$; that is, it is the bicategory of left $A$-, right $A$-bimodules. 
We have a commutative diagram
\begin{equation*}
\xymatrix{
\Delta^{\mathrm{op}}\ar[r]^-{\cong} \ar[d]_-{\widehat{A}}  & \Delta_{\bot, \top} \ar[r]^-{\mathrm{incl.}} & \Delta \ar[d]^-{A^*} \\
\tensor*[^A]{\CM}{^A} \ar[rr]_-{\mathrm{und.}} && \CM}
\end{equation*}
of pseudofunctors defining the augmented pseudosimplicial object $\widehat{A}$ of $\tensor*[^A]{\CM}{^A}$
which provides a (bicategorically) free resolution of $A$ acting on itself via
\begin{eqnarray}\label{p3}
P_3 = (A\ox A\ox A\xra{P\ox 1_A}A\ox A \xra{P}A) \ .
\end{eqnarray}
 Here is an indicative picture for low dimensions.
\begin{eqnarray}\label{resolution}
\xymatrix @R-3mm {\dots \
A^{\ox 4} \ar@<4ex>[rr]^{1_A\ox 1_A\ox P}  \ar@<0ex>[rr]^{1_A\ox P\ox 1_A} \ar@<-4ex>[rr]^{P\ox 1_A\ox 1_A} &&
A^{\ox 3} \ar@<3ex>[rr]^{1_{A}\ox P}  \ar@<-3ex>[rr]^{P\ox 1_A} && 
A^{\ox 2} \ar@<0.0ex>[rr]^{P} \ar@<0ex>[ll]_{1_A\ox J\ox 1_A}&& A  
}\end{eqnarray}
For any object $U\in \CM$, we obtain a pseudosimplicial object $U\ox \widehat{A}$ of $\CM$. 

Suppose however that $U$ is equipped with the structure $\underline{U} = (U,u,\zeta)$ of an object of the centre $\mathrm{Z}\CM$. Then we have a family of pseudo-equivalences 
\begin{eqnarray}\label{gradequiv}
u_A\ox 1_{A^{\ox (n-1)}} : U\ox A^{\ox n} \to A\ox U\ox A^{\ox (n-1)}
\end{eqnarray}
for $n\ge 2$.    
Dispensing with the augmentation, we can transport the pseudosimplicial structure on $U\ox \widehat{A}$
via the graded pseudo-equivalences \eqref{gradequiv} to obtain a pseudosimplicial object
\begin{eqnarray}\label{transported}
\xymatrix @R-3mm {\dots \
A\ox U \ox A^{\ox 3} \ar@<4ex>[rrr]^{1_{A\ox U\ox A}\ox P}  \ar@<0ex>[rrr]^{1_{A\ox U}\ox P\ox 1_A} \ar@<-4ex>[rrr]^{(P\ox 1_{U\ox A\ox A})\circ (1_A\ox u_A\ox 1_{A\ox A})} &&&
A\ox U \ox A^{\ox 2} \ar@<3ex>[rrr]^{1_{A\ox U}\ox P}  \ar@<-3ex>[rrr]^{(P\ox 1_{U\ox A})\circ (1_A\ox u_A\ox 1_A)} &&& 
A\ox U \ox A  \ar@<0ex>[lll]_{1_{A\ox U}\ox J\ox 1_A}  
}\end{eqnarray} 
of free objects in $\tensor*[^A]{\CM}{^A}$. 
Having the free structure on $A\ox X \ox A$ means, for any $M\in \tensor*[^A]{\CM}{^A}$, we have the equivalence of categories
\begin{eqnarray}\label{twosidedfree}
\tensor*[^A]{\CM}{^A}(A\ox X\ox A,M) \simeq \CM(X,M)
\end{eqnarray}
 given by composing with $J\ox 1_X\ox J : X\to A\ox X\ox A$. 
 Apply $\tensor*[^A]{\CM}{^A}(-,A)$ to the transported diagram \eqref{transported}, where
 $A$ has action as in \eqref{p3}. Now use the equivalences \eqref{twosidedfree} to
 obtain a pseudocosimplicial category \eqref{pscosimpcat} in which the functors $\partial_i$ and $\sigma_j$ are defined by \eqref{partsig} for $h : U\to A$ and $k : U\ox A\to A$. 
 \begin{eqnarray}\label{pscosimpcat}
\xymatrix @R-3mm {
\CM(U,A) \ar@<2.5ex>[rr]^{\partial_0\phantom{aa}}  \ar@<-2.5ex>[rr]^{\partial_1\phantom{aa}} &&
\CM(U\ox A,A) \ar@<0ex>[ll]_{\sigma_0\phantom{aa}}  \ar@<3ex>[rr]^{\partial_0\phantom{aaa}}
 \ar@<0ex>[rr]^{\partial_1\phantom{aaa}}   \ar@<-3ex>[rr]^{\partial_2\phantom{aaa}} && 
\CM(U\ox A\ox A,A) \ \dots  
}\end{eqnarray}   
\begin{eqnarray}\label{partsig}
\begin{aligned}
& \partial_0(h) = P\circ (h\ox 1_A)  \ & \  \partial_1(h) = P\circ (1_A\ox h)\circ u_A \phantom{aaaaaa}  
\\ & \sigma_0(k) = k\circ (1_U\ox J)  & \partial_0(k) = P\circ (k\ox 1_{A})  \phantom{aaaaaaaaaa}\\
& \  \partial_1(k) = k\circ (1_U\ox P) \  
& \phantom{aaaaaa} \partial_2(k) = P\circ (1_A\ox k)\circ (u_A\ox 1_A) 
\end{aligned} 
\end{eqnarray}

The category $\mathrm{CP}{\CM}(\underline{U},A)$ of {\em centre pieces} is the (strong) descent category \cite{31, 80}
for the pseudocosimplicial category \eqref{pscosimpcat}. 
This generalizes slightly the definition of \cite{81} since $\CM$ need not be braided
and we only require $\underline{U}\in \mathrm{Z}\CM$. An object of $\mathrm{CP}{\CM}(\underline{U},A)$
is a morphism $h : U\to A$ equipped with an invertible 2-morphism
\begin{equation}\label{cp2mor}
\begin{aligned}
\xymatrix{
U\ox A \ar[d]_{h\ox 1_A}^(0.8){\phantom{aaaaaaa}}="1" \ar[rr]^{u_A}  && A\ox U\ar[d]^{1_A\ox h}_(0.8){\phantom{aaaaaaa}}="2" \ar@{=>}"1";"2"^-{\gamma}
\\
A\ox A \ar[rd]_-{P} && A\ox A \ar[ld]^-{P} 
\\
& A 
}
\end{aligned}
\end{equation}
subject to a descent condition as in Section 2 of \cite{81}. The equation
\begin{eqnarray}\label{cpunitcondn}
\lgroup h \xRa{\rho h} P(h\ox j) \xRa{\gamma} P(J\ox h)\xRa{\lambda h}h \rgroup = 1_h
\end{eqnarray}
is a consequence.
The lax descent category for the pseudocosimplicial category \eqref{pscosimpcat} is
denoted by $\mathrm{CP}_{\ell}{\CM}(\underline{U},A)$; the objects are called {\em lax centre pieces} and differ only from centre pieces as described above in that the 2-morphism $\gamma$ in \eqref{cp2mor} is not required to be invertible and condition
\eqref{cpunitcondn} must be added as part of the descent condition; compare Section 4, Figure 1 of \cite{Ignacio}. 

The centre piece construction extends canonically to a pseudofunctor
\begin{eqnarray}\label{cp_psfun}
\mathrm{CP}{\CM}(-,A) : (\mathrm{Z}\CM)^{\mathrm{op}}\to \mathrm{Cat} \ .
\end{eqnarray}
A birepresenting object $\underline{\mathrm{Z}}A\in \mathrm{Z}\CM$ for \eqref{cp_psfun} is called the {\em full monoidal centre}
of the monoidale $A\in \CM$: this means we have equivalences
\begin{eqnarray*}
\mathrm{Z}\CM(\underline{U},\underline{\mathrm{Z}}A) \simeq \mathrm{CP}{\CM}(\underline{U},A)  
\end{eqnarray*}
pseudonatural in $\underline{U} \in \mathrm{Z}\CM$. 
So we have a universal centre piece $\mathrm{z}_A : \mathrm{Z}A\to A$
with a 2-morphism as in \eqref{cp2mor} with $\underline{U} = \underline{\mathrm{Z}}A$. 
The proof in \cite{81} of 
Proposition 2.1
(or the alternative in Remark 3.3 suggested by Stephen Lack) carries over to show that $\underline{\mathrm{Z}}A$ is a braided monoidale in $\mathrm{Z}\CM$ and
$\mathrm{z}_A : \mathrm{Z}A\to A$ is strong monoidal in $\CM$.
Of course, the {\em full monoidal lax centre}
of the monoidale $A\in \CM$ is defined by a psudonatural family of equivalences
\begin{eqnarray*}
\mathrm{Z}\CM(\underline{U},\underline{\mathrm{Z}}_{\ell}A) \simeq \mathrm{CP}_{\ell} {\CM}(\underline{U},A)  
\end{eqnarray*}

If $\CM$ is a monoidal category (that is, has only identity 2-morphisms) then 
$\underline{\mathrm{Z}}A$ is the full centre of the monoid $A$ in the sense of Davydov \cite{Davydov2010}.

Suppose now that $\CM$ is a braided monoidal bicategory and $A$ is a monoidale in $\CM$. Then each object $U$ of $\CM$ becomes an object $\underline{U}$ via the pseudofunctor \eqref{MinZM} and we write $\mathrm{CP}{\CM}(U,A)$ rather than
$\mathrm{CP}{\CM}(\underline{U},A)$ for the category of centre pieces. 
Restriction of \eqref{cp_psfun} along \eqref{MinZM} provides a pseudofunctor  
\begin{eqnarray}\label{restrictedcp}
\mathrm{CP}{\CM}(-,A) : \CM^{\mathrm{op}}\to \mathrm{Cat} \ .
\end{eqnarray}
In agreement with \cite{81}, a birepresenting object $\mathrm{Z}A\in \CM$ for \eqref{restrictedcp} is called the {\em monoidal centre}
of the monoidale $A\in \CM$: this means we have equivalences
\begin{eqnarray*}
\CM(U,\mathrm{Z}A) \simeq \mathrm{CP}{\CM}(U,A)  
\end{eqnarray*}
pseudonatural in $U \in \CM$. This $\mathrm{Z}A$ is braided in $\CM$ and we
have a universal strong monoidal centre piece $\mathrm{z} : \mathrm{Z}A\to A$
in $\CM$.   

If $\CM$ is a braided closed monoidal bicategory, as mentioned in \cite{81}, the monoidal centre of a monoidale $A\in\CM$ is the codescent object for the
pseudosimplicial object 
\begin{eqnarray}\label{homversion}
\xymatrix @R-3mm {
[A\ox A, A] \ \ar@<3ex>[rr]^{d_0}  \ar@<0ex>[rr]^{d_1} \ar@<-3ex>[rr]^{d_2} &&
\ [A,A] \ \ar@<2ex>[rr]^{d_0}  \ar@<-2ex>[rr]^{d_1} && 
\ A \ .  \ar@<0ex>[ll]_{s_0}  
}\end{eqnarray} 
This follows by replacing each $\CM(U\ox A^{\ox n},A)$ in \eqref{transported} by the
pseudonaturally equivalent $\CM(U, [A^{\ox n},A])$ and applying the bicategorical Yoneda lemma of \cite{14}. 

Of course, the centre of a monoidal category $\CV$ in the sense of \cite{38} is
the monoidal centre $\mathrm{Z}\CV$ of the monoidale $\CV$ in
the braided closed cartesian monoidal bicategory $\mathrm{Cat}$.

\section{Centre pieces in $\mathrm{Mod}^{\mathrm{op}}$}\label{CpiModop}

Lax centres in $\CV\text{-}\mathrm{Mod}^{\mathrm{op}}$ were analysed in Section 7.2 of \cite{Ignacio}. We provide a brief reminder as needed for our purpose.

\begin{lemma}\label{natsatcol} 
Suppose that $S, T : \CA\to \CX$ and $K :\CD \to \CA$ and  are functors. Suppose that
every object $A$ of $\CA$ is a colimit of some functor which factors through $K$ and
that these colimits are preserved by $S$. Then restriction along $K$ provides a
bijection
$$[\CA,\CX](S,T) \cong [\CD,\CX](SK,TK) \ .$$
Moreover, if $T$ also preserves the colimits in question then the bijection restricts to the invertible natural transformations on both sides.  
\end{lemma}
\begin{proof} 
For $D : \CK\to \CD$ and $\theta \in [\CA,\CX](S,T)$, we see that 
 $\theta_{\mathrm{colim}KD} : S\mathrm{colim}KD\to T\mathrm{colim}KD$
 is induced on the colimit by the components of $\theta_{KD}$. 
 This allows the unique reconstruction of $\theta$ from any natural transformation 
 purporting to be its restriction. Moreover, any morphism induced on colimits by
 an invertible natural transformation is invertible. 
\end{proof}   

\begin{proposition}\label{cpModCat} 
For any monoidale $A$ in $\mathrm{Mod}^{\mathrm{op}}$, the isomorphism of categories $\mathrm{Mod}^{\mathrm{op}}(U,A)\cong [U^{\mathrm{op}},[A,\mathrm{Set}]]$ induces
isomorphisms of categories
\begin{eqnarray*}
\mathrm{CP}{\mathrm{Mod}^{\mathrm{op}}}(U,A) \cong \mathrm{CP}{\mathrm{Cat}}(U^{\mathrm{op}}, [A,\mathrm{Set}]) \ , \ \mathrm{CP}_{\ell}{\mathrm{Mod}^{\mathrm{op}}}(U,A) \cong \mathrm{CP}_{\ell}{\mathrm{Cat}}(U^{\mathrm{op}}, [A,\mathrm{Set}]) 
\end{eqnarray*}
where $[A,\mathrm{Set}]$ has the convolution monoidal structure for the promonoidal category  $A$. 
The three isomorphisms are pseudonatural in $U\in \mathrm{Cat}^{\mathrm{op}}$ where, on the left-hand sides, $U$ is inserted into $\mathrm{Mod}^{\mathrm{op}}$ via the top 
pseudofunctor $(-)^*$ of \eqref{starmodules} and, on the right-hand sides, $U$ is inserted into $\mathrm{Cat}$ via the left 2-functor $(-)^{\mathrm{op}}$ of \eqref{starmodules}.  
\end{proposition} 
\begin{proof}
We identify the functor $P : A^{\mathrm{op}}\x A^{\mathrm{op}} \to [A,\mathrm{Set}]$ with the module $P : A\to A\x A$ involved in the monoidale structure on $A$.  
A centre piece $(h,\widehat{\gamma}) : U^{\mathrm{op}}\to [A,\mathrm{Set}]$ in $\mathrm{Cat}$
consists of a functor $h : U^{\mathrm{op}}\to [A,\mathrm{Set}]$, which clearly
identifies with a module $h : A\to U$, and a natural isomorphism 
\begin{equation*} 
\begin{aligned}
\xymatrix{
U^{\mathrm{op}}\x [A,\mathrm{Set}] \ar[d]_{h\ox 1_{[A,\mathrm{Set}]}}^(0.8){\phantom{aaaaaaaaaaaaaaa}}="1" \ar[rr]^{c_{U^{\mathrm{op}},[A,\mathrm{Set}]}}  && [A,\mathrm{Set}]\x U^{\mathrm{op}}\ar[d]^{1_{[A,\mathrm{Set}]}\ox h}_(0.8){\phantom{aaaaaaaaaaaaaaa}}="2" \ar@{=>}"1";"2"^-{\widehat{\gamma}}
\\
[A,\mathrm{Set}]\x [A,\mathrm{Set}] \ar[rd]_-{\widehat{P}} && [A,\mathrm{Set}]\x [A,\mathrm{Set}] \ar[ld]^-{\widehat{P}} 
\\
& [A,\mathrm{Set}] 
}
\end{aligned}
\end{equation*}
satisfying the descent data condition. Here $\widehat{P}$ is colimit preserving in each
variable and has its restriction along $\mathrm{y}_A \x \mathrm{y}_A$, where 
$\mathrm{y}_A : A^{\mathrm{op}}\to [A,\mathrm{Set}]$ is the Yoneda embedding, isomorphic to the functor $P$.
As $\mathrm{y}_A$ is dense and both domain and codomain functors of $\widehat{\gamma}$ preserve colimits in the second variable,
we can apply Lemma~\ref{natsatcol} with respect to restriction along $U^{\mathrm{op}}\x \mathrm{y}_A$. This gives the bijection between
(invertible) natural transformations $\widehat{\gamma}$ and (invertible) 2-cells 
$\gamma : P\circ (h\x 1_A)\Rightarrow P\circ (1_A\x h)\circ c_{U,A}$.
Furthermore, the descent conditions appropriately correspond.      
\end{proof}     

\begin{remark}
The proof of Proposition~\ref{cpModCat} implicitly uses the fact that, if a pseudonatural transformation $f : M\to N$ between pseudocosimplicial categories $M$ and $N$ is such that $f_0$ is an equivalence, $f_1$ is fully faithful, and $f_2$ is faithful,
then $f$ induces an equivalence between the descent categories of $M$ and $N$. 
\end{remark}

Lax colimits exist in the bicategory $\mathrm{Mod}$ (for example, see \cite{85})
and this includes lax codescent objects. 
So each promonoidal category $A$ has a lax monoidal centre in $\mathrm{Mod}^{\mathrm{op}}$. 

\begin{corollary}\label{double} 
The lax centre $\mathrm{Z}_{\ell}A$ of a monoidale $A$ in $\mathrm{Mod}^{\mathrm{op}}$ satisfies a lax-braided monoidal equivalence
\begin{eqnarray*}
[\mathrm{Z}_{\ell}A, \mathrm{Set}] \simeq \mathrm{Z}_{\ell}[A,\mathrm{Set}] \ .
\end{eqnarray*}
If the centre $\mathrm{Z}A$ of a monoidale $A$ exists in $\mathrm{Mod}^{\mathrm{op}}$ then there is a braided monoidal equivalence
\begin{eqnarray*}
[\mathrm{Z}A, \mathrm{Set}] \simeq \mathrm{Z}[A,\mathrm{Set}] \ .
\end{eqnarray*}
\end{corollary} 
\begin{proof} The composite equivalence
\begin{align*}
[U^{\mathrm{op}}, \mathrm{Z}_{\ell}[A,\mathrm{Set}]] & \cong
\mathrm{CP}_{\ell}{\mathrm{Cat}}(U^{\mathrm{op}}, [A,\mathrm{Set}])\cong
\mathrm{CP}_{\ell}{\mathrm{Mod}}^{\mathrm{op}}(U, A) \\
& \simeq 
{\mathrm{Mod}}^{\mathrm{op}}(U, \mathrm{Z}_{\ell}A)
\cong [U^{\mathrm{op}}, [\mathrm{Z}_{\ell}A,\mathrm{Set}]]  
\end{align*}
gives the result on application of the bicategorical Yoneda Lemma.
 \end{proof}
\begin{remark}\label{lax=strong}
 If $A$ is an autonomous monoidal category then the lax monoidal centre is equivalent to the monoidal centre in $\mathrm{Mod}^{\mathrm{op}}$ (see \cite{89, Ignacio, 99}). 
\end{remark}

\section{The groupoid of automorphisms in a groupoid}\label{Tgoaiag}

For any groupoid $\CG$ and any autonomous monoidal category $\CV$, the category
$[\CG,\CV]$ of representations of $\CG$ in $\CV$, with the pointwise tensor product,
is also autonomous.

The groupoid of automorphisms in a groupoid $\CG$ will be denoted by $\CG^{\mathrm{aut}}$. 
The objects are pairs $(p,a)$ where $a : p\to p$ in $\CG$. 
A morphism $f : (p,a)\to (q,b)$ is a morphism $f\ :p\to q$ in $\CG$ such that $f a = b f$. 
We have a discrete fibration $\mathrm{q}_{\CG} : \CG^{\mathrm{aut}} \to \CG$ taking
 $(p,a)\xra{f}(q,b)$ to $p\xra{f}q$, with a section $\mathrm{i}_{\CG} : \CG\to \CG^{\mathrm{aut}}$ taking $p\xra{g}q$ 
 to $(p,1_p)\xra{g} (q,1_q)$.  

If $\mathrm{Aut}_{\CG} : \CG\to \mathrm{Set}$ denotes the functor taking the morphism $p\xra{g}q$ in $\CG$ to the conjugation function
$\CG(p,p)\xra{\CG(g^{-1},g)} \CG(q,q), \ a\mapsto \tensor*[^g]{a}{}$, then there is a standard equivalence of categories
\begin{eqnarray}\label{slice_equiv}
[\CG,\mathrm{Set}]/\mathrm{Aut}_{\CG} \simeq [\CG^{\mathrm{aut}}, \mathrm{Set}]
\end{eqnarray}
because $\CG^{\mathrm{aut}}$ is the category of elements of $\mathrm{Aut}_{\CG}$: a natural transformation $\phi : X\to \mathrm{Aut}_{\CG}$ corresponds to 
the functor $\phi^* : \CG^{\mathrm{aut}}\to \mathrm{Set}$ which takes the object $(p,a)$ to the fibre of 
$\phi_p : Xp \to \CG(p,p)$ over $a$ and takes $(p,a)\xra{f}(q,b)$ to the restriction
of $Xp\xra{Xf}Xq$ to those fibres.        

There is a monoid structure on $\mathrm{Aut}_{\CG}$ in the cartesian monoidal category $[\CG,\mathrm{Set}]$
given by componentwise composition in $\CG$. Consequently, there is a monoidal structure
on the left hand side of \eqref{slice_equiv} whose tensor product takes cartesian product of the
morphisms over $\mathrm{Aut}_{\CG}$ followed by the monoid multiplication;
also, there is a braiding as pointed out by Freyd and Yetter \cite{FY}.
This monoidal structure is closed (on both sides) and so transports to a promonoidal structure on
$\CG^{\mathrm{aut}}$: recall from Example 9 in Section 7 of \cite{60} or the end of Section 4 of \cite{87} that the promonoidal structure is defined by
\begin{eqnarray}\label{promononG^aut}
P((p,a),(q,b); (r,c)) = \{ p\xra{u}r\xla{v}q :  \tensor*[^u]{a}{} \tensor*[^v]{b}{} = c \} & \text{and} 
& J(r,c) = \left\{
\begin{array}{ll}
1 & \text{if  } c = 1_r \\
\varnothing & \text{if  } c \neq 1_r 
\end{array} \right. 
\end{eqnarray}
and that there is a braiding 
\begin{eqnarray*}
\gamma_{a,b;c} : P((p,a),(q,b); (r,c)) \xra{\cong} P((q,b),(p,a);(r,c)) \ , & (u,v)\mapsto (\tensor*[^u]{a}{} v,u) \ . 
\end{eqnarray*}
It also has a twist
$$\tau_a = a : (p,a) \to (p,a) \ ;$$
compare Section 2 of \cite{Street1998}.
The reader is invited to check the commutativity of \eqref{twistcondition} which is the main twist condition.
\begin{equation}\label{twistcondition}
\begin{aligned}
\xymatrix{
P((p,a),(q,b); (r,c)) \ar[rr]^-{\gamma_{a,b;c}} \ar[d]_-{P(1,1;\tau_c)} && P((q,b),(p,a);(r,c)) \ar[d]^-{P(\tau_a,\tau_b;1)} \\
P((p,a),(q,b); (r,c)) && P((q,b),(p,a);(r,c)) \ar[ll]_-{\gamma_{b,a;c}} }
\end{aligned}
\end{equation}
Furthermore, $\CG^{\mathrm{aut}}$ is a $*$-autonomous promonoidal category in the sense of \cite{79}:
we have the natural isomorphisms
\begin{eqnarray*}
P((p,a),(q,b); (r,c^{-1})) \xra{\cong} P((q,b),(r,c); (p,a^{-1}))  \ , & \ 
(u,v) \leftrightsquigarrow (u^{-1}v,u^{-1}) \ ;
\end{eqnarray*}
\begin{eqnarray*}
P((p,a),(q,b); (r,c^{-1})) \xra{\cong} P((q,b^{-1}),(p,a^{-1}); (r,c))^{\mathrm{op}} \ ,  &  \
(u,v) \leftrightsquigarrow (v,u) \ .
\end{eqnarray*}

From \cite{81, 87, 89}, we extract:

\begin{proposition}\label{centresforG}
The convolution braided monoidal category $[\CG^{\mathrm{aut}},\mathrm{Set}]$ 
is braided monoidal equivalent to the monoidal centre $\mathrm{Z}[\CG,\mathrm{Set}]$ of the 
cartesian monoidal category $[\CG,\mathrm{Set}]$. 
\end{proposition}
The monoidal equivalence of Proposition~\ref{centresforG} is the composite of the equivalence $[\CG,\mathrm{Set}]/\mathrm{Aut}_{\CG} \simeq \mathrm{Z}[\CG,\mathrm{Set}]$ taking $X \xra{\phi} \mathrm{Aut}_{\CG}$ to $(X, X\times - \xra{u} -\times X)$,
where $u_{Y, p}(x,y) = (Y(\phi_p(x))y,x)$, and the standard equivalence \eqref{slice_equiv}.   
We note that the equivalence becomes balanced on transport of the convolution twist to the monoidal centre.
\begin{remark}\label{on:centresforG}
The universal centre piece for the monoidal centre of $[\CG,\mathrm{Set}]$ composed
with the equivalence of Proposition~\ref{centresforG} is isomorphic
to the composite 
\begin{eqnarray*}
[\CG^{\mathrm{aut}},\mathrm{Set}]\simeq [\CG,\mathrm{Set}]/\mathrm{Aut}_{\CG} \xra{\mathrm{dom}} [\CG,\mathrm{Set}] 
\end{eqnarray*}
which is left Kan extension $\mathrm{Lan}_{\rm{q}}$ along the discrete fibration $\mathrm{q}_{\CG} : \CG^{\mathrm{aut}}\to \CG$.
On objects, it takes the functor $S : \CG^{\mathrm{aut}}\to \mathrm{Set}$ to the
functor $\hat{S} = \mathrm{Lan}_{\rm{q}}S : \CG \to \mathrm{Set}$ defined by $\hat{S}p = \sum_{a\in \CG(p,p)}S(p,a)$ with $\hat{S}(p\xra{g}q)$ taking $x\in S(p,a)$ to
$(Sg)x\in S(q,\tensor*[^g]{a}{})$. The centre piece structure $\gamma$ on $\hat{S}$ has components
\begin{eqnarray*}
 \gamma_{S,Y, p} : \sum_{a\in \CG(p,p)}S(p,a)\times Yp\lra Yp\times \sum_{a\in \CG(p,p)}S(p,a) \\ 
 (x\in S(p,a),y\in Yp)\mapsto (Y(a)y\in Yp, x\in S(p,a)) \ . 
\end{eqnarray*}
\end{remark}

\begin{proposition}\label{Gautascentre}
The braided monoidale $\CG^{\mathrm{aut}}$ is the monoidal centre of the monoidale $\CG$ 
in $\mathrm{Mod}^{\mathrm{op}}$.
The braided monoidal bicategory $\mathrm{Ps}(\CG^{\mathrm{aut}},\mathrm{Mod}^{\mathrm{op}})$ is
the monoidal centre of the autonomous monoidal bicategory $\mathrm{Ps}(\CG,\mathrm{Mod}^{\mathrm{op}})$.
\end{proposition}
\begin{proof}
The first sentence follows from Corollary~\ref{double}, Remark~\ref{lax=strong} and Proposition~\ref{centresforG}. The second sentence is a mildly higher dimensional version of Theorem 8.6 in \cite{89} and Remark~\ref{on:centresforG}. Notice that, for any set 
$\Lambda$, the conservative left biadjoint of the diagonal pseudofunctor 
$\mathrm{Diag} : \mathrm{Mod}^{\mathrm{op}} \to \mathrm{Ps}(\Lambda,\mathrm{Mod}^{\mathrm{op}})$ 
is also a right biadjoint since bicategorical coproducts in $\mathrm{Mod}$ are also bicategorical products \cite{85}; so 
$\Lambda$ is decomposing for $\mathrm{Mod}^{\mathrm{op}}$ in the language of \cite{89}. 
(We could replace $\CG$ by any category in which all endomorphisms are invertible). 
\end{proof}

\begin{remark}\label{on:Gautascentre}
The explicit biequivalence
\begin{eqnarray}\label{centrebiequiv}
\mathrm{Z}\mathrm{Ps}(\CG,\mathrm{Mod}^{\mathrm{op}}) \sim \mathrm{Ps}(\CG^{\mathrm{aut}},\mathrm{Mod}^{\mathrm{op}})
\end{eqnarray}
was well prepared for in Example~\ref{centreofptwise}. 
Take an object $(F,u, \zeta)$ of the centre of the form \eqref{deltavariant} involving $\delta_p : Fp\to \CG(p,p)$.
Define $\check{F} : \CG^{\mathrm{aut}}\to \mathrm{Mod}^{\mathrm{op}}$ on objects by taking the fibre 
$\check{F}(p,a) = \delta_p^{-1}(a)$ of $\delta_p$ over $a\in \CG(p,p)$. On the morphism $g : (p,a)\to (q,b)$ it is
the restriction of the functor $Fg$ to the fibres which makes sense using the dinaturality \eqref{dinatdelta}.   
The inverse biequivalence takes $S : \CG^{\mathrm{aut}}\to\mathrm{Mod}^{\mathrm{op}}$ to 
$(\hat{S},\delta)$ where $\hat{S}p= \sum_{a\in\CG(p,p)}{S(p,a)}$ and $\delta_p$ picks off the index.
This means that each component $u_{T,p} : Tp\times \hat{S}p\to \hat{S}p\times Tp$ 
of the centre structure is represented on the direct sums by a diagonal matrix with entries
\begin{eqnarray}\label{diagonalentries}
u_{T,p,a} : Tp\times S(p,a)\to S(p,a)\times Tp , \phantom{aaaaa} (\tau, \sigma) \mapsto (\sigma,(Ta)\tau) \ .
\end{eqnarray}
\end{remark}

\begin{example}\label{ex2}
Return now to our fibration $\pi : \CH\to \CG$ of Example~\ref{ex1}. 
Since fibrations in $\mathrm{Cat}$ are
preserved by 2-functors of the form $[\CD,-]$, we have the fibration $\pi^{\mathrm{aut}} : \CH^{\mathrm{aut}}\to \CG^{\mathrm{aut}}$. 
The corresponding pseudofunctor $\mathbb{H}^{\mathrm{aut}} : \CG^{\mathrm{aut}}\to \mathrm{Mod}^{\mathrm{op}}$ is defined as follows. 
The category $\mathbb{H}^{\mathrm{aut}} (p,a)$ has objects those $(s,x)\in \CH^{\mathrm{aut}}$ with $\pi(x)=a$;
morphisms $(s,x) \xra{k} (s_1,x_1)$ are those in $\CH^{\mathrm{aut}}$ with $\pi(k) = 1_p$. 
For $(p,a) \xra{f} (q,b)$ in $\CG^{\mathrm{aut}}$, the module $\mathbb{H}^{\mathrm{aut}}f : \mathbb{H}^{\mathrm{aut}} (q,b)\to \mathbb{H}^{\mathrm{aut}} (p,a)$ is defined by
\begin{eqnarray*}
(\mathbb{H}^{\mathrm{aut}} f)((s,x),(t,y)) = \{h\in \CH^{\mathrm{aut}}((s,x),(t,y)) : \pi(h) = f \} \ .
\end{eqnarray*}

The cleavage $\sigma$ for $\pi$ also gives a cleavage for $\pi^{\mathrm{aut}}$.
The corresponding pseudofunctor $\mathbb{H}^{\mathrm{aut}} : \CG^{\mathrm{aut}}\to \mathrm{Cat}^{\mathrm{op}}$ has the same value on objects as in the last paragraph. 
For $(p,a) \xra{f} (q,b)$ in $\CG^{\mathrm{aut}}$, the functor $\mathbb{H}^{\mathrm{aut}}f : \mathbb{H}^{\mathrm{aut}} (q,b)\to \mathbb{H}^{\mathrm{aut}} (p,a)$
takes $(t,y)\xra{k} (t_1,y_1)$ in $\mathbb{H}^{\mathrm{aut}} (q,b)$ to $(f^*(t),y^{\sigma(f)}) \xra{f^*(k)} (f^*(t_1),y_1^{\sigma(f)})$.
The invertible composition constraint $(\mathbb{H}^{\mathrm{aut}} f)(\mathbb{H}^{\mathrm{aut}} g)\Ra \mathbb{H}^{\mathrm{aut}}(gf)$ has component at $(u,z)\in \mathbb{H}^{\mathrm{aut}} (r,c)$ equal to
$\sigma(gf)^{-1}\sigma(g)\sigma(f) : (f^*g^*(u), z^{\sigma(g)\sigma(f)})\to ((gf)^*(u),z^{\sigma(gf)})$.
It is worth remembering that the functor $\mathbb{H}^{\mathrm{aut}}f$ is an equivalence with pseudo-inverse $\mathbb{H}^{\mathrm{aut}}f^{-1}$. 

We have the pseudofunctor $\widehat{\mathbb{H}^{\mathrm{aut}}} : \CG\to \mathrm{Cat}^{\mathrm{op}}$, 
supporting the centre structure, corresponding
to $\mathbb{H}^{\mathrm{aut}}$ under the biequivalence \eqref{centrebiequiv}. 
The commutative triangle
\begin{equation*}
\xymatrix{
\CH^{\mathrm{aut}} \ar[rd]_{ }\ar[rr]^{\mathrm{q}_{\CH}}   && \CH \ar[ld]^{\pi} \\
& \CG  &
}
\end{equation*}
of fibrations induces functors $\mathrm{q}_{\mathbb{H} p} : \widehat{\mathbb{H}^{\mathrm{aut}}}p\to \mathbb{H}p$
between the fibres. The right adjoint modules $\mathrm{q}_{\mathbb{H} p}^* : \mathbb{H}p \to \widehat{\mathbb{H}^{\mathrm{aut}}}p$ as morphisms in $\mathrm{Mod}^{\mathrm{op}}$ are the components of a morphism
\begin{eqnarray}\label{z for bbH}
\mathrm{z}_{\mathbb{H}} : \widehat{\mathbb{H}^{\mathrm{aut}}}\to \mathbb{H}
\end{eqnarray} 
in $\mathrm{Ps}(\CG,\mathrm{Mod}^{\mathrm{op}})$.   
\end{example} 

 \section{Monoidales in convolution bicategories}\label{Micb}

One virtue of the promonoidal groupoid $\CG^{\mathrm{aut}}$ over the monoid $\mathrm{Aut}_{\CG}$ is
that we can obtain convolution balanced monoidal structures on functors from $\CG^{\mathrm{aut}}$, not only into
$\mathrm{Set}$ but, into any nice enough monoidal category; or even on pseudofunctors from 
$\CG^{\mathrm{aut}}$ into any nice enough monoidal bicategory.     

Let $\CK$ be a monoidal bicategory with coproducts preserved by horizontal composition in each variable.
The tensor product will be denoted by $-\ox - : \CK\times \CK\to \CK$ with unit object $\CI$. 
Think of $\CG^{\mathrm{aut}}$ as a bicategory with only identity 2-cells.
We will make explicit the convolution monoidal structure on the bicategory 
$\mathrm{Ps}(\CG^{\mathrm{aut}},\CK)$.

Take $S, T\in \mathrm{Ps}(\CG^{\mathrm{aut}},\CK)$. 
Put\footnote{The coend here is in the pseudo-sense appropriate to bicategories.}
\begin{eqnarray*}
(S\star T)(r,c) = \sum_{ab =c}{S(r,a)\ox T(r,b)} 
\ \bigl(\simeq \int^{(p,a),(q,b)}_{\mathrm{ps}}{P((p,a),(q,b); (r,c))\cdot S(p,a)\ox T(q,b)} \ \bigr)
\end{eqnarray*}
and, for $(r,c)\xra{f} (r_1,c_1)$ in $\CG^{\mathrm{aut}}$, define $(S\star T)f$ by commutativity in
\begin{equation}
\begin{aligned}
\xymatrix{
S(r,a)\ox T(r,b)\ar[rr]^-{\mathrm{in}_{a,b}} \ar[d]_-{Sf\ox Tf} &&  \sum_{ab =c}{S(r,a)\ox T(r,b)} \ar[d]^-{(S\star T)f} \\
S(r_1, \tensor*[^f]{a}{})\ox T(r_1, \tensor*[^f]{b}{}) \ar[rr]_-{\mathrm{in}_{\tensor*[^f]{a}{},\tensor*[^f]{b}{}}} && \sum_{a_1b_1 =c_1}{S(r_1,a_1)\ox T(r_1, b_1)} \ .}
\end{aligned}
\end{equation}
This defines the tensor product $S\star T$ for a monoidal structure on
$\mathrm{Ps}(\CG^{\mathrm{aut}},\CK)$ with unit $\mathbb{J} : \CG^{\mathrm{aut}}\to \CK$ defined by
\begin{eqnarray}
\begin{aligned}
\mathbb{J}(r,c) = \left\{
\begin{array}{ll}
\CI & \text{if  } c = 1_r \\
0 & \text{if  } c \neq 1_r 
\end{array} \right.
\end{aligned}
\end{eqnarray}
 which becomes functorial on noting that, for $c\xra{f} c'$, if $c=1_r$
 then $c'=1_r$. 
 
 We can now contemplate monoidales $M$ in $\mathrm{Ps}(\CG^{\mathrm{aut}},\CK)$. 
 Such a monoidale consists of a pseudofunctor $M : \CG^{\mathrm{aut}}\to \CK$ equipped with morphisms
 \begin{eqnarray*}
I_r : \CI\to M(r,1_r) \ \text{  and  } \ \square_{a,b} : M(r,a)\ox M(r,b)\to M(r,ab)
\end{eqnarray*}
 in $\CK$ and invertible 2-cells
 \begin{equation}\label{assoconstraints}
 \begin{aligned}
\xymatrix{
M(r,a)\ox M(r,b)\ox M(r,c) \ar[d]_{\square_{a,b}\ox 1}^(0.5){\phantom{aaaaaaaaa}}="1" \ar[rr]^{1\ox \square_{b,c}}  && M(r,a)\ox M(r,bc) \ar[d]^{\square_{a,bc}}_(0.5){\phantom{aaaaaaaaa}}="2" \ar@{=>}"1";"2"^-{\alpha_{a,b,c}}_{\cong}
\\
M(r,ab)\ox M(r,c) \ar[rr]_-{\square_{ab,c}} && M(r,abc) 
}
\end{aligned}
\end{equation}  
\begin{equation}\label{unitconstraints}
 \begin{aligned}
   \xymatrix{ & & M(r,a) \ar[dd]^{1} \ar[lld]_{I_r\ox 1} \ar[rrd]^{1\ox I_r}\\
     M(r,1_r)\ox M(r,a)   \ar[rrd]_{\square_{1,a}} & & \rtwocell[0.4]{ll}{\lambda_a \cong} & & \rtwocell{ll}{\rho_a \cong} M(r,a)\ox M(r,1_r) \ar[lld]^{\square_{a,1}}\\
     & & M(r,a)}
 \end{aligned}
\end{equation}
all subject to pseudonaturality 
\begin{equation}\label{pseudonatofpromonoidalM}
\begin{aligned}
\xymatrix{\ar @{} [drr] | {\stackrel{\cong} \Longrightarrow}
M(r,a)\ox M(r,b)\ar[rr]^-{\square_{a,b}} \ar[d]_-{Mf\ox Mf} &&  M(r,ab) \ar[d]^-{Mf} \\
M(r,\tensor*[^f]{a}{})\ox M(r,\tensor*[^f]{b}{}) \ar[rr]_-{\square_{\tensor*[^f]{a}{},\tensor*[^f]{b}{}}} &  & M(r,\tensor*[^f]{a}{}\tensor*[^f]{b}{}) \ ,} 
\end{aligned}
\end{equation}
modificationality and coherence conditions. 

\begin{example}\label{ex3}
When we consider the monoidal bicategory $\mathrm{Ps}(\CG,\mathrm{Mod}^{\mathrm{op}})$,
it is with the pointwise monoidal structure, which is autonomous since $\mathrm{Mod}^{\mathrm{op}}$ is and $\CG$ is a groupoid. 
 The $\mathbb{H}$ of Example~\ref{ex1} is a monoidale in $\mathrm{Ps}(\CG,\mathrm{Mod}^{\mathrm{op}})$. 
 The monoidal structure is provided by the modules 
 $$\square_p :  \mathbb{H}p\to \mathbb{H}p \times \mathbb{H}p \ \text{ and } \ I_p :   \mathbb{H}p \to \mathbf{1}$$
 defined by  
 \begin{eqnarray*}
\square_p(s, t; u) = \mathbb{H}p(s,u)\times \mathbb{H}p(t,u)\ , & \ \square_p(x, y; z) = \mathbb{H}p(x,z)\times \mathbb{H}p(y,z) \ , & \ I_p = !_* 
\end{eqnarray*}
The unit and associativity constraints are much as for the promonoidal structure \eqref{promonoidalG} on $\CG$.
The pseudonaturality structure on the modules $\square(s, t; u)$ is provided by the Yoneda Lemma isomorphisms
\begin{align*}
((\mathbb{H}a \times \mathbb{H}a)\circ \square_q)(s_1,s_2;t) & =  \int^{t_1,t_2}{\square_q(t_1, t_2; t)\times \mathbb{H}q(s_1,a^*(t_1))\times \mathbb{H}q(s_2,a^*(t_2))} \\
&  = \int^{t_1,t_2}{\mathbb{H}q(t_1,t)\times \mathbb{H}q(t_2,t)\times \mathbb{H}p(s_1,a^*(t_1))\times \mathbb{H}p(s_2,a^*(t_2))} \\
& \cong \ \mathbb{H}p(s_1,a^*(t))\times \mathbb{H}p(s_2,a^*(t)) = \square_p(s_1, s_2; a^*(t)) \\
& \cong \int^s \mathbb{H}p(s,a^*(t))\times \square_p(s_1, s_2; s) = (\square_p\circ \mathbb{H}a)(s_1,s_2;t) \ . 
\end{align*}
\end{example}
\begin{example}\label{ex4}
 The $\mathbb{H}^{\mathrm{aut}}$ of Example~\ref{ex2} is a monoidale in $\mathrm{Ps}(\CG^{\mathrm{aut}},\mathrm{Mod}^{\mathrm{op}})$. 
 The monoidal structure is provided by the modules 
 $$\square_{a,b} :  \mathbb{H}^{\mathrm{aut}}(r,ab)\to \mathbb{H}^{\mathrm{aut}} (r,a) \times \mathbb{H}^{\mathrm{aut}} (r,b) \ \text{ and } \ I_r :   \mathbb{H}^{\mathrm{aut}} (r,1_r) \to \mathbf{1}$$
 defined by  
 \begin{eqnarray*}
 \square_{a,b}((s,x), (t,y); (u,z)) = \{(m,n)\in \mathbb{H}r(s,u)\times \mathbb{H}r(t,u) : \tensor*[^m]{x}{}\tensor*[^n]{y}{} =z \}
 \end{eqnarray*}
  \begin{eqnarray*}
I_r(u,z) = \left\{
\begin{array}{ll}
1 & \text{if  } z = 1_u \\
\varnothing & \text{if  } z \neq 1_u 
\end{array} \right. \ . 
\end{eqnarray*}
The unit and associativity constraints are much as for the promonoidal structure \eqref{promononG^aut} on $\CG^{\mathrm{aut}}$.
The pseudonaturality structure on the component modules $\square_{a,b}((s,x), (t,y); (u,z))$ is provided by the Yoneda Lemma isomorphisms
\begin{align*}
((\mathbb{H}^{\mathrm{aut}}f \times \mathbb{H}^{\mathrm{aut}}f)\circ \square_{\tensor*[^f]{a}{},\tensor*[^f]{b}{}})((s,x),(t,y); (u,z)) & \cong  \square_{a,b}((s,x), (t,y); (f^*(u),z^{\sigma(f)})) \\
&  \cong  (\square_{a,b} \circ \mathbb{H}^{\mathrm{aut}}f)((s,x),(t,y); (u,z)) \ . 
\end{align*}
We leave the coend calculation as an exercise with the reminder that the three occurrences of the functor $\mathbb{H}^{\mathrm{aut}}f$
are values at the three morphisms $(p,a)\xra{f} (q,\tensor*[^f]{a}{})$, $(p,b)\xra{f} (q,\tensor*[^f]{b}{})$, and $(p,ab)\xra{f} (q,\tensor*[^f]{a}{}\tensor*[^f]{b}{})$ (as per \eqref{pseudonatofpromonoidalM}). 
 \end{example} 
\bigskip

\section{The Turaev-Virelizier structures}\label{tTVs}

The definitions in this section are those of \cite{TurVir2013}. For this section, fix a group $G$ regarded as a groupoid with
one object $o_G$ and with morphisms $o_G\xra{g} o_G$ the elements $g$ of the group.
 
\begin{definition}A {\em $G$-graded category over $k$} is a $k$-linear monoidal category $\CC$, with finite direct sums, endowed with a system of pairwise disjoint full $k$-linear subcategories $\CC_{a}, a\in G$, with finite direct sums, such that
\begin{itemize}
\item[(a)] each object $X\in \CC$ splits as a direct sum $\oplus_a{X_a}$ where $X_a\in \CC_a$ and
$a$ runs over a finite subset of $G$;
\item[(b)] if $X\in \CC_a$ and $Y\in \CC_b$ then $X\ox Y\in \CC_{ab}$;
\item[(c)] if $X\in \CC_a$ and $Y\in \CC_b$ with $a\ne b$ then $\CC(X,Y)=0$;
\item[(d)] the tensor unit $I$ of $\CC$ is in $\CC_1$.
\end{itemize} 
\end{definition}
\bigskip

Turaev-Virelizier call an object $X$ of a $G$-graded category $\CC$ {\em homogeneous} when there exists a 
(necessarily unique) $a\in G$ such that $X\in \CC_a$; this $a$ is denoted by $|X|$.

They write $\bar{G}$ for the discrete monoidal category of elements of 
$G$ with the multiplication as tensor product.
They write $\mathrm{Aut}(\CC)$ for the monoidal category of monoidal endo-equivalences of
the monoidal $k$-linear category $\CC$ and monoidal natural isomorphisms; 
the tensor product is composition of functors.
\bigskip

\begin{definition} A $G$-crossed category $\CC$ is a 
$G$-graded category over $k$ equipped with a strong monoidal functor $\phi : \bar{G}\to \mathrm{Aut}(\CC)$
such that $\phi_a(\CC_b) \subseteq \CC_{a^{-1}ba}$ for all $a,b\in G$.
\end{definition}
\bigskip

\begin{proposition}\label{convmonprop}
Let $\CV$ be the monoidal category of modules over a fixed commutative ring. Then monoidales in $\mathrm{Ps}(G^{\mathrm{aut}},\CV\text{-}\mathrm{Cat})$ 
are equivalent to the $G$-crossed categories of \cite{TurVir2013}.
\end{proposition}
\begin{proof}
Take a monoidale $M$ in $\mathrm{Ps}(G^{\mathrm{aut}},\CV\text{-}\mathrm{Cat})$. 
Let $\CC$ be the $\CV$-category obtained by taking the completion, with respect to finite direct sums, of the coproduct $\CC_{\mathrm{hom}}= \sum_a{Ma}$ in $\CV\text{-}\mathrm{Cat}$. 
Then $\CC$ is a $G$-graded category with $\phi_f = Mf^{-1}$.

Conversely, take a $G$-graded category $\CC$ and define $Ma$ to be the full sub-$\CV$-category $\CC_a$ of $\CC$ consisting of the objects homogeneous over $a\in G$.
Then $M$ is a monoidale in $\mathrm{Ps}(G^{\mathrm{aut}},\CV\text{-}\mathrm{Cat})$.      
\end{proof}

\begin{corollary}
With $\CV$ as in Proposition~\ref{convmonprop}, any monoidale in $\mathrm{Ps}(G^{\mathrm{aut}},\mathrm{Mod}^{\mathrm{op}})$ delivers a $G$-crossed category on application of the pseudofunctor \eqref{[-,V]}.
\end{corollary}

\begin{definition} A {\em $G$-braided category $\CC$} is a 
$G$-crossed category $(\CC, \phi)$ equipped with a natural family of isomorphisms
\begin{eqnarray*}
\gamma_{X,Y} : X\ox Y\lra Y\ox \phi_{|Y|}(X) \ ,
\end{eqnarray*}
for $X, Y\in \CC$ and $Y$ homogeneous, subject to \textbf{three} axioms. 
\end{definition}
In the next section, we will see how this fits into our theory of braided monoidales.

\section{Internal homs, biduals and braidings}\label{Ihbab}

If $\CK$ is left closed, it is straightforward to see that so too is $\mathrm{Ps}(\CG^{\mathrm{aut}},\CK)$:
 \begin{eqnarray}\label{closed}
[T, U](p,a) = \prod_b{[T(p,b),U(p,ab)]} \ . 
\end{eqnarray}
\begin{proposition}
Suppose in $\CK$ that direct sums indexed by the endohomsets of $\CG$ exist and that each $T(q,b)$ has a left bidual $T(q,b)^{\vee}$. Then 
$T$ has a left bidual
 \begin{eqnarray}\label{bidual}
T^{\vee}(p,a) = [T, \mathbb{J}](p,a) = T(p,a^{-1})^{\vee} 
\end{eqnarray}
in $\mathrm{Ps}(\CG^{\mathrm{aut}},\CK)$.
\end{proposition}
\begin{proof}
Taking \eqref{bidual} as the definition of $T^{\vee}$, we need to prove that the 
canonical morphism $S\star T^{\vee}\lra [T,S]$
is an equivalence for all $S$. The component of this canonical pseudonatural
transformation at $(r,c)$ is the composite
 \begin{eqnarray*}
\sum_{ab=c}S(r,a)\ox T(r,b^{-1})^{\vee} \simeq \sum_{d}S(r,cd)\ox T(r,d)^{\vee} \\
\simeq \sum_{d}{[T(r,d), S(r,cd)]} \xra{\mathrm{canon.}} \prod_{d}{[T(r,d), S(r,cd)]} 
\end{eqnarray*}
in which the arrow is an equivalence because of our assumption about direct sums.
\end{proof}
\begin{corollary}
All biduals exist
in $\mathrm{Ps}(\CG^{\mathrm{aut}},\mathrm{Mod}^{\mathrm{op}})$; 
that is, the monoidal bicategory is autonomous (also called ``compact'' or ``rigid'').\end{corollary}

 If $\CK$ is equipped with a braiding $\gamma_{X,Y} : X\ox Y\to Y\ox X$ then we obtain a braiding $\gamma_{S,T} : S\star T\to T\star S$ on $\mathrm{Ps}(\CG^{\mathrm{aut}},\CK)$
as defined by the commutative pentagon \eqref{convbraid}.
\begin{eqnarray}\label{convbraid}
\begin{aligned}
\xymatrix{
S(r,a)\ox T(r,b) \ar[r]^-{\gamma_{S(r,a),T(r,b)}} \ar[d]_-{\mathrm{in}_{a,b}} & T(r,b)\ox S(r,a) \ar[r]^-{ Ta\ox 1}  & T(r,\tensor*[^a]{b}{})\ox S(r,a) \ar[d]^-{\mathrm{in}_{\tensor*[^a]{b}{},a}} \\
\sum_{ab =c}{S(r,a)\ox T(r,b)} \ar[rr]_-{\gamma_{S,T, (r,c)}} & & \sum_{b'a' =c}{T(r,b')\ox S(r,a')} }
\end{aligned}
\end{eqnarray}
\bigskip

If $\CK$ is balanced then so too is $\mathrm{Ps}(\CG^{\mathrm{aut}},\CK)$ with twist 
$$\theta_{S, (p,a)} = \bigl(S(p,a) \xra{Sa} S(p,a) \xra{\theta_{S(p,a)}} S(p,a)\bigr) \ .$$
Recall that, if $\CK$ is symmetric, we choose its twist to be the identity.   
 
We already have the example $\CG^{\mathrm{aut}}$ of a $*$-autonomous
balanced monoidale in $\mathrm{Mod}^{\mathrm{op}}$. 

\begin{proposition}\label{tortile}
The monoidal bicategory $\mathrm{Ps}(\CG^{\mathrm{aut}},\mathrm{Mod}^{\mathrm{op}})$ is tortile 
(in the sense of \cite{ShumPhD}).
\end{proposition}

Moreover, with $\mathrm{Ps}(\CG^{\mathrm{aut}},\CK)$ a braided monoidal bicategory, 
according to \cite{81}, we can contemplate
monoidal centres $\mathrm{Z}M$ for monoidales $M$ therein.
Since the centre is a limit, it is formed pointwise in $\CK$.
From \cite{81}, we know that $\mathrm{Z}M$ is a braided monoidale in $\mathrm{Ps}(\CG^{\mathrm{aut}},\CK)$. 
\bigskip   

\begin{example}\label{ex5}
 The monoidale $\mathbb{H}^{\mathrm{aut}}$ of Example~\ref{ex4} is balanced in $\mathrm{Ps}(\CG^{\mathrm{aut}},\mathrm{Mod}^{\mathrm{op}})$. 
The braiding is
\begin{eqnarray*}
\gamma_{x,y;z} : \square_{a,b}((s,x), (t,y); (u,z)) \xra{\cong} \square_{b,a}((t,y), (s,x); (u,z)) \ , & (m,n)\mapsto (\tensor*[^m]{x}{} n,m) \ .
\end{eqnarray*}
The twist $\tau : \mathbb{H}^{\mathrm{aut}} \to \mathbb{H}^{\mathrm{aut}}$ is given by $\tau_a = \mathbb{H}^{\mathrm{aut}} a : \mathbb{H}^{\mathrm{aut}} (p,a) \to \mathbb{H}^{\mathrm{aut}} (p,a)$.
 \end{example} 
\bigskip   
 
 For the remainder of this section, let us return to the context of Section~\ref{tTVs} where $G$ is a group and $\CK = \CV\text{-}\mathrm{Cat}$ where $\CV$ is a complete cocomplete closed symmetric monoidal category. 
For $S\in \mathrm{Ps}(G^{\mathrm{aut}},\CV\text{-}\mathrm{Cat})$,
$f : a \to b$ in $G^{\mathrm{aut}}$ and $A\in S$, we put $fA = (Sf)A \in Sb$.   

Let $M$ be a monoidale in $\mathrm{Ps}(G^{\mathrm{aut}},\CV\text{-}\mathrm{Cat})$. 
The tensor product consists of $\CV$-functors $\square_{a,b} : Ma\ox Mb\to M(ab)$. The unit is an object $I$ of $M1$.
The associativity constraint consists of a $\CV$-natural family
\begin{eqnarray}\label{assoconstrcomponents}
\begin{aligned}
\alpha_{a,b,c} : (A\square_{a,b} B)\square_{ab,c} C\lra A\square_{a,bc}(B\square_{b,c} C) \ .
\end{aligned}
\end{eqnarray}
A braiding for $M$ consists of a $\CV$-natural family 
\begin{eqnarray}\label{braidingcomponents}
\begin{aligned}
\gamma_{a,b} : A\square_{a,b} B\lra aB \square_{\tensor*[^a]{b}{},a} A \ .
\end{aligned}
\end{eqnarray} 

\begin{proposition}\label{convbraidprop}
In the setting of Proposition~\ref{convmonprop}, the braided monoidales in \\
$\mathrm{Ps}(G^{\mathrm{aut}},\CV\text{-}\mathrm{Cat})$ 
are equivalent to the $G$-braided categories of \cite{TurVir2013}.
\end{proposition}
\bigskip

According to Section 3 of \cite{81}, since pseudolimits limits are formed pointwise, 
the monoidal centre $\mathrm{Z}M$ of a monoidale $M$ in $\mathrm{Ps}(G^{\mathrm{aut}},\CV\text{-}\mathrm{Cat})$ is constructed as follows.
The $\CV$-category $(\mathrm{Z}M)a$ has objects pairs $(A,\upsilon)$
where $A$ is an object of $Ma$ and $\upsilon$ is a {\em half $G$-braiding} for $A$
consisting of a $\CV$-natural family of isomorphisms 
\begin{eqnarray}\label{halfbraidingcomponents}
\begin{aligned}
\upsilon_{b} : A\square_{a,b} B\lra aB \square_{\tensor*[^a]{b}{},a} A 
\end{aligned}
\end{eqnarray}
such that $\upsilon_{1} : A\square_{a,1} I\lra I \square_{1,a} A$ transports the right unit constraint into
the left unit constraint and the following hexagon commutes.  
\begin{equation*}
\xymatrix{
& A\square_{a,bc} (B\square_{b,c} C) \ar[rr]^-{\upsilon}  & & (aB\square_{\tensor*[^a]{b}{},\tensor*[^a]{c}{}} aC)\square_{\tensor*[^a]{(bc)}{},a} A) \ar[d]^-{\alpha} & \\
(A\square_{a,b} B)\square_{ab,c} C \ar[ur]^-{\alpha} \ar[dr]_-{\upsilon\ox 1} & & & aB\square_{\tensor*[^a]{b}{},ac} (aC\square_{\tensor*[^a]{c}{},a} A) \\
 & (aB\square_{\tensor*[^a]{b}{},a} A)\square_{ab,c} C \ar[rr]_-{\alpha} & & aB\square_{\tensor*[^a]{b}{},ac} (A\square_{a,c} C) \ar[u]_-{1\ox \upsilon} &}
\end{equation*}
For $f : a \to \tensor*[^f]{a}{}$, we have $(\mathrm{Z}M)(f)(A,\upsilon) = (fA, \upsilon')$ where
$\upsilon'_b$ for $B\in Mb$  is the composite 
\begin{equation*}
fA\square_{\tensor*[^f]{a}{},b}B\cong f(A\square_{a,\tensor*[^f^{-1}]{b}{}}f^{-1}B)\xrightarrow{f\upsilon_{f^{-1}B}}f(af^{-1}B\square_{\tensor*[^a]{(f^{-1}b)}{},a}A)\cong \tensor*[^f]{a}{}B\square_{\tensor*[^{\tensor*[^f]{a}{}}]{b}{}, \tensor*[^f]{a}{}}fA  \ .
\end{equation*}

% \section{Biequivalents of $\mathrm{Ps}(\CG,\mathrm{Mod}^{\mathrm{op}})$ when $\CG$ is a groupoid}
% 
% Let $\mathrm{Cat}_{\mathrm{cc}}$ denote the full sub-2-category of $\mathrm{Cat}$
% consisting of the categories in which idempotents split. We have the locally fully faithful
% pseudofunctor
% \begin{eqnarray}\label{CatccinMod}
% \mathrm{K} : \mathrm{Cat}_{\mathrm{cc}} \lra \mathrm{Mod}^{\mathrm{op}} 
%\end{eqnarray}
%taking each category $\CA$ to its opposite $\CA^{\mathrm{op}}$ and each functor
%$F : \CA \to \CB$ to the module $F_* : \CB^{\mathrm{op}}\to \CA^{\mathrm{op}}$. 
%
%\begin{proposition}
%For any groupoid $\CG$ with blah blah blah, composition with the pseudofunctor \eqref{CatccinMod} defines a biequivalence of bicategories
%\begin{eqnarray*}
%\Theta : \mathrm{Ps}(\CG, \mathrm{Cat}_{\mathrm{cc}}) \xra{\sim} \mathrm{Ps}(\CG, \mathrm{Mod}^{\mathrm{op}}) \ .
%\end{eqnarray*}
%\end{proposition}
%\begin{proof}
%Take pseudofunctors $S, T : \CG \to \mathrm{Cat}_{\mathrm{cc}}$ and a pseudonatural transformation
%$\alpha : \mathrm{K}S \to \mathrm{K}T$.  
%\end{proof}

 \section{A full centre from a groupoid fibration}\label{Afcfagf}
 
 Davydov \cite{Davydov2010} defined the {\em full centre} of a monoid $A$ in a (not
 necessarily braided) monoidal category
 $\CV$ to be a commutative monoid $\underline{\mathrm{Z}}A$ in the (braided) monoidal centre $\mathrm{Z}\CV$ 
 of $\CV$ satisfying an appropriate universal property. The universal property is that of a terminal object in a category of elements of a set-valued functor $\mathrm{CP}\CV(\underline{-},A)$, and so amounts to saying that the set-valued functor is representable. It is pointed out in \cite{Street2012} that
 the pair $(\mathrm{Z}\CV, \underline{\mathrm{Z}}A)$ is the monoidal centre of the monoidale $(\CV, A)$
 in the monoidal bicategory of pointed categories.  
 
 In Section \ref{Roic}, we lifted the concept of full centre from the monoidal category level to the monoidal bicategory level.
 The full monoidal centre of a monoidale $A$ in a monoidal bicategory $\CM$ is a braided monoidale
 $\underline{\mathrm{Z}}M$ in the monoidal centre $\mathrm{Z}\CM$ of $\CM$ satisfying an appropriate universal property.    

\begin{proposition}
The full monoidal centre of the monoidale $\mathbb{H}$ in $\mathrm{Ps}(\CG,\mathrm{Mod}^{\mathrm{op}})$ (see Examples~\ref{ex1} and \ref{ex3}) is the braided monoidale $\mathbb{H}^{\mathrm{aut}}$ in
$\mathrm{Ps}(\CG^{\mathrm{aut}},\mathrm{Mod}^{\mathrm{op}})$ (see Examples~\ref{ex4} and \ref{ex5}).
The universal centre piece has underlying morphism $\mathrm{z}_{\mathbb{H}} : \widehat{\mathbb{H}^{\mathrm{aut}}}\to \mathbb{H}$ as in \eqref{z for bbH}. 
\end{proposition}
\begin{proof}
We must produce an equivalence 
\begin{eqnarray}\label{toprove(ii)}
\mathrm{Ps}(\CG^{\mathrm{aut}},\mathrm{Mod}^{\mathrm{op}})(S,\mathbb{H}^{\mathrm{aut}}) \simeq 
\mathrm{CP Ps}(\CG,\mathrm{Mod}^{\mathrm{op}})(\underline{\hat{S}},\mathbb{H})  
\end{eqnarray} 
where we can suppose that $S$ factors through $(-)_* : \mathrm{Cat}^{\mathrm{op}}\to\mathrm{Mod}^{\mathrm{op}}$
and where $\hat{S}$ is as in Remark~\ref{on:Gautascentre}.

 Recalling the centre structure \eqref{diagonalentries} on $\hat{S}$ and performing some coend calculations, 
 we see that the centre piece structure on a pseudonatural transformation $k : \hat{S}\to \mathbb{H}$ 
 amounts to isomorphisms 
\begin{eqnarray}\label{k_pa}
k_{p,a}(\sigma,s)\times \mathbb{H}p(s',a^*(s)) \cong \mathbb{H}p(s',s)\times k_{p,a}(\sigma,s)   
\end{eqnarray}
of the form $(\alpha, z)\mapsto (\kappa(\alpha,z), \alpha)$, where $k_{p,a}$ is the composite module 
$$\mathbb{H}p\xra{k_p}\sum_{a\in \CG(p.p)}S(p,a)\xra{\mathrm{in}^*_a}S(p,a) \ .$$   

Let $h : S\to \mathbb{H}^{\mathrm{aut}}$ be a pseudonatural transformation.
It has component modules $h_{(p,a)} : \mathbb{H}^{\mathrm{aut}}(p,a)\to S(p,a)$ for objects $(p,a) \in \CG^{\mathrm{aut}}$.  
Identify $h_{(p,a)}$ with the corresponding functor $\mathbb{H}^{\mathrm{aut}}(p,a)\to [S(p,a)^{\mathrm{op}},\mathrm{Set}]$ and let $\hat{h}_{p,a} : \mathbb{H}p\to [S(p,a)^{\mathrm{op}},\mathrm{Set}]$ be the left Kan extension of that functor
along the discrete fibration $\mathbb{H}^{\mathrm{aut}}(p,a)\to \mathbb{H}p$ taking $(s,x)$ to $s$. Explicitly,
\begin{eqnarray}\label{littlehhat}
\hat{h}_{p,a}(\sigma,s) = \sum_{{\scriptsize \begin{matrix}  x\in \CH(s,s) \\ \pi(x) = a \end{matrix}}}{h_{(p,a)}(\sigma,(s,x))} 
\end{eqnarray}
 for each $\sigma\in S(p,a)$. 
For \eqref{k_pa} with $k=\hat{h}$, we take $\kappa(\alpha,z) = x^{-1}\sigma(s)z$ where $x\in \CH(s,s)$ is the index of the summand of 
\eqref{littlehhat} which contains $\alpha$. This define \eqref{toprove(ii)} on objects.

Now take a morphism $\omega : h\Ra h' : S\to \mathbb{H}^{\mathrm{aut}}$.
We obtain natural transformations $\omega_{(p,a)} : h_{(p,a)}\Ra h'_{(p,a)} : \mathbb{H}^{\mathrm{aut}}(p,a)\to [S(p,a)^{\mathrm{op}},\mathrm{Set}]$ which induce natural transformations $\hat{\omega}_{p,a} : \hat{h}_{p,a}\Ra \hat{h'}_{p,a}$ on the Kan extensions leading to a modification $\hat{\omega} : \hat{h} \Ra \hat{h'}$ compatible with the centre
piece structures. 

Thus we have a family of functors 
$$(\widehat{-})_S : \mathrm{Ps}(\CG^{\mathrm{aut}},\mathrm{Mod}^{\mathrm{op}})(S,\mathbb{H}^{\mathrm{aut}}) \to 
\mathrm{CP Ps}(\CG,\mathrm{Mod}^{\mathrm{op}})(\underline{\hat{S}},\mathbb{H})$$ 
pseudonatural in $S$.
By the bicategorical Yoneda Lemma, the family is induced by the $\hat{h}$ obtained by setting $S = \mathbb{H}^{\mathrm{aut}}$ and taking $h$ to be the identity module. Then the corresponding functor
$h_{(p,a)}$ is the Yoneda embedding. The module $\hat{h}_p : \mathbb{H}p\to \widehat{\mathbb{H}^{\mathrm{aut}}}p$
is the right adjoint $\mathrm{q}_{\mathbb{H} p}^*$ mentioned in the line before \eqref{z for bbH}. 
Therefore $\hat{h}$ is the morphism 
$\mathrm{z}_{\mathbb{H}} :  \widehat{\mathbb{H}^{\mathrm{aut}}} \to \mathbb{H}$ equipped with its centre piece structure.

Now take a pseudonatural transformation $k : \hat{S}\to \mathbb{H}$
equipped with a centre piece structure determined as in \eqref{k_pa} by functions $\kappa : k_{p,a}(\sigma,s)\times \mathbb{H}p(s',a^*(s)) \to \mathbb{H}p(s',s)$. Define modules $\check{k}_{(p,a)} : \mathbb{H}^{\mathrm{aut}}(p,a)\to S(p,a)$ by  
\begin{eqnarray*}
\check{k}_{(p,a)}(\sigma,(s,x)) = \{\alpha \in k_{p,a}(\sigma,s) : \kappa(\alpha, 1_{a^*(s)}) = x^{-1}\sigma(s) \} 
\end{eqnarray*}
to obtain the components on objects of a pseudonatural transformation $\check{k} : S\to \mathbb{H}^{\mathrm{aut}}$.
We immediately see from \eqref{littlehhat} that $\hat{\check{k}}= k$.
Using naturality of \eqref{k_pa} in $s'$, we see that the $\kappa$ for $\check{\hat{h}}$ recovers that for $h$. 
The construction $k\mapsto \check{k}$ extends to a functor
$$(\widecheck{-})_S :  
\mathrm{CP Ps}(\CG,\mathrm{Mod}^{\mathrm{op}})(\underline{\hat{S}},\mathbb{H}) \to \mathrm{Ps}(\CG^{\mathrm{aut}},\mathrm{Mod}^{\mathrm{op}})(S,\mathbb{H}^{\mathrm{aut}})$$
which is an inverse equivalence for $(\widehat{-})_S$.         
\end{proof}
\bigskip
%
%\noindent \textbf{Exercises} 1. Repeat the whole story replacing $G$ by a Hopf algebra. 
%
%2. Repeat the whole story replacing $G$ by a Hopf monoidal comonoid.

%\bibliographystyle{acm}
%\bibliography{../../Notes/references_c.bib}

\begin{thebibliography}{1}

\bibitem{BN1996}
{\sc J.C. Baez and M. Neuchl,} 
\newblock Higher dimensional algebra, I. Braided monoidal 2-categories. 
\newblock {\em Advances in Math. 121} (1996) 196--244.

\bibitem{BarSch2011}
{\sc T. Barmeier and C. Schweigert,} 
\newblock A geometric construction for permutation equivariant categories from modular functors. 
\newblock {\em Transformation Groups 16} (2011) 287--337.

\bibitem{Ben1967} 
{\sc Jean B\'enabou,} 
\newblock Introduction to bicategories.
\newblock {\em Lecture Notes in Mathematics 47} (Springer-Verlag, 1967) 1--77.

\bibitem{BenLectures} 
{\sc Jean B\'enabou,} 
\newblock {\em Lectures at Oberwolfach (Germany) and other places} (1972--1995).

\bibitem{BenDistrib} 
{\sc Jean B\'enabou,} 
\newblock Les distributeurs.
\newblock {\em Universit\'e Catholique de Louvain, Inst. de Math. Pure et Appliqu\'es, Rapport 33} (1973).

\bibitem{Bruguieres2011}
{\sc A. Brugui{\'e}res, S. Lack and A. Virelizier,}
\newblock Hopf monads on monoidal categories.
\newblock {\em Advances in Mathematics 227(2)} (2011) 745--800.

\bibitem{Bunge1966} 
{\sc Marta Bunge,}
\newblock Categories of Set-Valued Functors. 
\newblock PhD Thesis (University of Pennsylvania 1966).

\bibitem{Bunge2011} 
{\sc Marta Bunge,}
\newblock Adjoint functors between diagrammatic categories. 
\newblock Note (Montr\'eal, 3 February 2011).

\bibitem{Cond} 
{\sc F. Conduch\'e,}  
\newblock Au sujet de l'existence d'adjoints \`a droite aux foncteurs ``image r\'eciproque'' dans la cat\'egorie des cat\'egories. 
\newblock {\em Compte Rendus Acad. Sci. Paris S\'er. A-B 275} (1972) A891--A894.

\bibitem{Chikhladze2010}
{\sc D. Chikhladze, S. Lack and R. Street,}
\newblock Hopf monoidal comonads.
\newblock {\em Theory Appl. Categ. 24 No. 19} (2010) 554--563.

\bibitem{Crans1998}
{\sc Sjoerd E. Crans,} 
\newblock Generalized centers of braided and sylleptic monoidal 2-categories. 
\newblock {\em Advances in Math. 136} (1998) 183--223.

\bibitem{Davydov2010}
{\sc Alexei Davydov,} 
\newblock Centre of an algebra.
{\em Advances in Mathematics 225} (2010) 319--348

\bibitem{Day1970}
{\sc Brian J. Day,} 
\newblock On closed categories of functors.  
\newblock{\em in ``Reports of the Midwest Category Seminar IV'', Lecture Notes in Mathematics 137} (Springer, Berlin, 1970) 1--38.

\bibitem{87} 
{\sc B.J. Day, E. Panchadcharam and R. Street,} 
\newblock Lax braidings and the lax centre. 
\newblock {\em in ``Hopf Algebras and Generalizations'', Contemporary Mathematics 441} (2007) 1--17.

\bibitem{60} 
{\sc B.J. Day and R. Street,}
\newblock Monoidal bicategories and Hopf algebroids. 
\newblock {\em Advances in Math. 129} (1997) 99--157.

\bibitem{79} 
{\sc B.J. Day and R. Street,}
\newblock Quantum categories, star autonomy, and quantum groupoids. 
\newblock {\em in "Galois Theory, Hopf Algebras, and Semiabelian Categories", Fields Institute Communications 43} (American Math. Soc. 2004) 187--226.

\bibitem{89} 
{\sc B.J. Day and R. Street,}
\newblock Centres of monoidal categories of functors.
\newblock {\em in ``Categories in Algebra, Geometry and Mathematical Physics'', Contemporary Mathematics 431} (2007) 187--202.
 
\bibitem{Drinfeld} 
{\sc V.G. Drinfel'd,} 
\newblock Quantum groups. 
\newblock {\em Proceedings of the International Congress of Mathematicians at Berkeley, California, U.S.A. 1986} (1987) 798--820.

\bibitem{EilKel1966} 
{\sc S. Eilenberg and G.M. Kelly,} 
\newblock Closed categories. 
\newblock {\em Proceedings of the Conference on Categorical Algebra (La Jolla, 1965)} (Springer-Verlag,1966) 421--562.

\bibitem{Freyd1966} 
{\sc Peter Freyd,} 
\newblock Abelian Categories. 
{\em Harper International Edition} (1966).

\bibitem{FY} 
{\sc P. Freyd and D. Yetter,} 
\newblock Braided compact closed categories with applications to low dimensional
topology. 
{\em Advances in Math. 77} (1989) 156--182.

\bibitem{FW} 
{\sc A. Fr\"ohlich and C. T. C. Wall,} 
\newblock Graded monoidal categories. 
{\em Compositio Mathematica 28(3)} (1974) 229--285.

\bibitem{Giraud64} 
{\sc Jean Giraud,} 
\newblock M\'ethode de la descente.
\newblock {\em Bull. Soc. Math. France M\'em. 2} (1964).

\bibitem{GPS} 
{\sc R. Gordon, A.J. Power and R. Street,} 
\newblock Coherence for tricategories.
\newblock {Memoirs of the American Mathematical Society 117(558)} (1995) vi+81 pp.

\bibitem{38} 
{\sc A. Joyal and R. Street,}
\newblock  Tortile Yang-Baxter operators in tensor categories. 
\newblock {\em J. Pure Appl. Algebra 71} (1991) 43--51.

\bibitem{Joyal1991}
{\sc A. Joyal and R. Street,}
\newblock An introduction to Tannaka duality and quantum groups.
\newblock In {\em Category Theory: Proceedings of the International Conference
  held in Como, Italy, July 22--28, 1990}, A.~Carboni, M.~C. Pedicchio, and
  G.~Rosolini, Eds., vol.~1488 of {\em Lecture Notes in Mathematics}. Springer
  Berlin Heidelberg, 1991, pp.~413--492.
  
\bibitem{Joyal1993}
{\sc A. Joyal and R. Street,}
\newblock Braided tensor categories.
\newblock {\em Advances in Math. 102} (1993) 20--78.

\bibitem{Kelly1974}
{\sc M. Kelly and R. Street,} 
\newblock Review of the elements of 2-categories.
\newblock {\em Lecture Notes in Math. 420} (1974) 75--103.

\bibitem{KirPri2008}
{\sc A. Kirillov and T. Prince,} 
\newblock On G-modular functor.
\newblock arXiv:0807.0939

\bibitem{KTZ}
{\sc  L. Kong, Y. Tian and S. Zhou,}
\newblock The center of monoidal 2-categories in 3+1D Dijkgraaf-Witten Theory.
\newblock {\em Advances in Math. 360} (2020) 106--928.

\bibitem{LawMetric} 
{\sc  F. W. Lawvere,} 
\newblock Metric spaces, generalized logic and closed categories. 
\newblock {\em Reprints in Theory and Applications of Categories 1} (2002) pp.1--37.

\bibitem{Ignacio}
{\sc  Ignacio L\'opez Franco} 
\newblock Formal Hopf algebra theory. II. Lax centres. 
\newblock {\em J. Pure Appl. Algebra 213(11)} (2009) 2038--2054.

\bibitem{99}
{\sc  I. L\'opez Franco, R.H. Street and R.J. Wood,} 
\newblock Duals invert.
\newblock {\em Applied Categorical Structures 19(1)} (2011) 321--362. 

\bibitem{MacLPar} 
{\sc S. Mac Lane and R. Par\'e,} 
\newblock Coherence for bicategories and indexed categories. 
\newblock {\em Journal of Pure and Applied Algebra 37} (1985) 59--80.

\bibitem{MNS}
{\sc J. Maier, T. Nikolaus and C. Schweigert,} 
\newblock Equivariant modular categories via Dijkgraaf-Witten theory. 
\newblock arXiv:1103.2963

\bibitem{Majid}
{\sc Shahn Majid,} 
\newblock Representations, duals and quantum doubles of monoidal categories. 
\newblock {\em Rendiconti del Circolo Matematico di Palermo. Serie II. Supplemento. No. 26} (1991) 197--206.

\bibitem{Paddy}
{\sc Paddy McCrudden,}
\newblock Balanced coalgebroids. 
\newblock {\em Theory and Applications of Categories 7(6)} (2000) 71--147.

\bibitem{BMitchell} 
{\sc Barry Mitchell,}
\newblock Rings with several objects. 
\newblock {\em Advances in Math. 8} (1972) 1--161.

\bibitem{Mueger}
{\sc Michael M\"uger,} 
\newblock From subfactors to categories and topology. II. The quantum double of tensor categories and subfactors. 
\newblock {\em J. Pure Appl. Algebra 180} (2003) 159--219.

\bibitem{ShumPhD}
{\sc Mei Chee Shum,}
\newblock Tortile tensor categories.
\newblock {\em Journal of Pure and Applied Algebra 93} (1994) 57-110.

\bibitem{Street1972a}
{\sc Ross Street,}
\newblock The formal theory of monads.
\newblock {\em Journal of Pure and Applied Algebra 2(2)} (1972) 149--168.

\bibitem{14}
{\sc Ross Street,}
\newblock Fibrations in bicategories.
\newblock {\em Cahiers de topologie et g\'eom\'etrie diff\'erentielle 21} (1980) 111--160.

\bibitem{85}
{\sc Ross Street,}
\newblock Cauchy characterization of enriched categories. 
\newblock {\em Reprints in Theory and Applications of Categories 4} (2004) 1--16.

\bibitem{31} 
{\sc Ross Street,}
\newblock Correction to ``Fibrations in bicategories''. 
\newblock {\em Cahiers de topologie et g\'eom\'etrie diff\'erentielle cat\'egoriques 28} (1987) 53--56.

\bibitem{Street1998}
{\sc Ross Street,}
\newblock The quantum double and related constructions. 
\newblock {\em J. Pure Appl. Algebra 132} (1998) 195--206

\bibitem{Street2001}
{\sc Ross Street,}
\newblock Powerful functors.
\newblock \url{http://web.science.mq.edu.au/~street/Pow.fun.pdf}

\bibitem{80}
{\sc Ross Street,}
\newblock Categorical and combinatorial aspects of descent theory. 
\newblock {\em Applied Categorical Structures 12} (2004) 537--576. 

\bibitem{81}
{\sc Ross Street,}
\newblock The monoidal centre as a limit. 
\newblock {\em Theory and Applications of Categories 13(13)} (2004) 184--190.

\bibitem{Street2007}
{\sc Ross Street,}
\newblock {\em Quantum Groups: A Path to Current Algebra}.
\newblock Australian Mathematical Society Lecture Series 19. Cambridge
  University Press, 2007.
  
\bibitem{Street2012}
{\sc Ross Street,}
\newblock Monoidal categories in, and linking, geometry and algebra. 
\newblock {\em Bulletin of the Belgian Mathematical Society -- Simon Stevin 19(5)} (2012) 769--821.  

\bibitem{Tur2000}
{\sc Vladimir Turaev,} 
\newblock Homotopy field theory in dimension 3 and crossed group-categories.
\newblock arxiv:math.GT/0005291
  
  \bibitem{TurVir2013}
{\sc V. Turaev and A. Virelizier,}
\newblock On the graded center of graded categories. 
\newblock {\em Journal of Pure and Applied Algebra 217} (2013) 1895--1941.

\end{thebibliography}

\end{document}